%%rettet 1.12.16 og 18.12.16, og 24.3.17, og 30.4.17, 25.2.18
\input  amstex
\input amsppt.sty
\magnification1200
\vsize=23.5truecm
\hsize=16.5truecm
%\vcorrection{-10truemm}
\NoBlackBoxes

\def\supp{\operatorname{supp}}

\def\crp{\overline{\Bbb R}_+}
\def\crm{\overline{\Bbb R}_-}
\def\crpm{\overline{\Bbb R}_\pm}

\def\rnp{{\Bbb R}^n_+}

\def\rnpm{\Bbb R^n_\pm}
\def\crnp{\overline{\Bbb R}^n_+}
\def\crnm{\overline{\Bbb R}^n_-}
\def\crnpm{\overline{\Bbb R}^n_\pm}
\def\comega{\overline\Omega }

\def\Rn{\Bbb R^n}
\def\ang#1{\langle {#1} \rangle}

\def\Op{\operatorname{Op}}
\def\Pfrac{\tsize\frac1{\raise 1pt\hbox{$\scriptstyle p$}}}
\def\pfrac{\frac1{\raise 1pt\hbox{$\scriptscriptstyle p$}}}
\def\Pfracc#1{\tsize\frac{#1}{\raise 1pt\hbox{$\scriptstyle p$}}}
\def\pfracc#1{\frac{#1}{\raise 1pt\hbox{$\scriptscriptstyle p$}}}

\def\simto{\overset\sim\to\rightarrow}

\def\Zfrac{\tsize\frac1{\raise 1pt\hbox{$\scriptstyle z$}}}
\def\zfrac{\frac1{\raise 1pt\hbox{$\scriptscriptstyle z$}}}

\def\rp{ \Bbb R_+}

\def\OP{\operatorname{OP}}

\def\R{\Bbb R}

\def\ol{\overline}
\def\SD{\Cal S}
\def\E{\Cal E}
\def\F{\Cal F}
\def\D{\Cal D}

\document
\topmatter
\title
Green's formula and a Dirichlet-to-Neumann operator for
fractional-order pseudodifferential operators 
\endtitle
\author Gerd Grubb \endauthor
\affil
{Department of Mathematical Sciences, Copenhagen University,
Universitetsparken 5, DK-2100 Copenhagen, Denmark.
E-mail {\tt grubb\@math.ku.dk}}\endaffil
\rightheadtext{Green's formula}
\abstract
The paper treats boundary value problems for the fractional Laplacian
$(-\Delta )^a$, $a>0$, and  more generally for
classical pseudodifferential operators ($\psi $do's) $P$ of order $2a$ with
even symbol, applied to functions on a smooth subset $\Omega $ of
${\Bbb R}^n$. There are several meaningful local boundary conditions,
such as the Dirichlet and Neumann conditions $\gamma _k^{a-1}u=\varphi
$, $k=0,1$, where $\gamma _k^{a-1}u=c_k\partial_n^k(u/d^{a-1})|_{\partial\Omega }$, $d(x)=\operatorname{dist}(x,\partial\Omega )$. We show
a new Green's formula $$(Pu,v)_\Omega -(u,P^*v)_\Omega =(s_0\gamma _1^{a-1}u+B\gamma
_0^{a-1}u,\gamma _0^{a-1}v)_{\partial\Omega }-(s_0\gamma _0^{a-1}u,\gamma _1^{a-1}v)_{\partial\Omega },$$
where $B$ is a first-order $\psi $do on $\partial\Omega $ depending on
the first two terms in the symbol of $P$. 

Moreover, we show in the elliptic case how the Poisson-like solution operator $K_D$ for the
nonhomogeneous Dirichlet problem is constructed from $P^+$ in the
factorization $P\sim P^-P^+$ obtained in earlier work. The
Dirichlet-to-Neumann operator $S_{DN}=\gamma _1^{a-1}K_D$ is derived
from this as a first-order $\psi $do on
$\partial\Omega $, with an explicit formula for the symbol. This leads to a characterization of those operators
$P$ for which the Neumann problem is Fredholm solvable.  
\endabstract
\endtopmatter

\subhead 1. Introduction \endsubhead

The fractional Laplacian $(-\Delta )^a$ on ${\Bbb R}^n$, $0<a<1$, and its
 boundary value problems on subsets $\Omega \subset{\Bbb R}^n$, have received
 much attention recently, as a generalization from the ordinary
 Laplacian $-\Delta $, with useful applications to
 probability, finance, differential geometry and mathematical
 physics. %[refs]
$(-\Delta )^a$ can be described as a pseudodifferential operator ($\psi $do) or as
a singular integral operator:
$$
\aligned
(-\Delta )^au&=\operatorname{OP}(|\xi |^{2a})u=
\Cal F^{-1}(|\xi |^{2a}\hat u(\xi ))\\
&=c_{n,a}PV\int_{{\Bbb R}^n}\frac{u(x)-u(y)}{|y|^{n+2a}}\,dy.
\endaligned\tag1.1
$$
The discussion we give in this paper works for classical $\psi
$do's $P$ of order $2a$, for {\it any} $a>0$, with symbol $p\sim\sum_{j\in{\Bbb
    N}_0}p_j(x,\xi )$ being {\it even}:
$$
p_j(x,-\xi )=(-1)^jp_j(x,\xi ),\text{ all }j.\tag1.2
$$
E.g.\ $P=A(x,D)^a$,
where $A(x,D)$ is a second-order strongly elliptic differential operator. 

One of the difficulties with these operators is that they are {\it
  nonlocal}, in contrast to differential operators.
This is problematic when one wants to study them on  subsets $\Omega \subset\Bbb
R^n$. 

It has been known for many years that one can define a Dirichlet
realization $P_D$ by a variational construction: Define the
sesquilinear form $p_0(u,v)$ by
$$
p_0(u,v)=\int_\Omega Pu\,\bar v\, dx,\quad u,v\in C_0^\infty (\Omega ),\tag1.3
$$
completed to a form on $\dot H^a(\overline\Omega )=\{u\in H^a(\Bbb R^n)\mid
 \operatorname{supp}u\subset\overline\Omega \}$. It is
coercive when $P$ is strongly elliptic, and
then induces an operator $P_D$ in $L_2(\Omega )$ acting like $r^+P$ with domain
$$
D(P_D)=\{u\in \dot H^a(\overline\Omega )\mid r^+Pu\in L_2(\Omega
)\},\quad r^+Pu=(Pu)|_\Omega ;\tag1.4
$$
a Fredholm operator when $\Omega $ is bounded.
% Here $r^+$ indicates restriction to $\Omega $.  
$P_D$ represents the (fractional restricted) homogeneous Dirichlet problem
$$
r^+Pu= f \text{ in }\Omega ,\quad \operatorname{supp}u\subset
\overline\Omega .\tag1.5
$$

Vishik and Eskin treated such operators in the 1960's (see e.g.\
 \cite{E81})
by pseudodifferential factorization methods; one consequence is that 
$D(P_D)= \dot
H^{2a }(\overline\Omega )$ when $a<\frac12$, and  $D(P_D)\subset\dot
H^{a+\frac12-\varepsilon  }(\overline\Omega )$ when $a\ge
\frac12$. 

More recently, real integral operator methods have been applied. There
is a technique by Caffarelli
and Silvestre \cite{CS07} to view $(-\Delta )^a$ on $\Bbb R^n$ as the Dirichlet-to-Neumann
operator for a degenerate elliptic differential boundary value
problem  on $\Bbb R^n\times \Bbb R_+$, which allows local techniques.
For subsets $\Omega \subset {\Bbb R}^n$,  
methods from potential theory and probability have particularly been used, leading to results in
H\"older spaces, under low smoothness assumptions. Also functional
analysis methods enter. 
Let us mention some of the studies through the times:  
Blumenthal and Getoor \cite{BG59}, Landkof \cite{L72}, Hoh and Jacob
\cite{HJ96},  Kulczycki
\cite{K97}, Chen and Song \cite{CS98}, Jakubowski \cite{J02},
Silvestre \cite{S07},
Caffarelli and Silvestre \cite{CS09}, Musina and Nazarov \cite{MN14}, Frank and Geisinger
\cite{FG16}, 
Ros-Oton and Serra \cite{RS14a, RS14b, RS16}, Abatangelo \cite{A15}, Felsinger, Kassmann and Voigt
\cite{FKV15}, Bonforte, Sire
and Vazquez \cite{BSV15}, Servadei and Valdinoci \cite{SV14}, Ros-Oton
\cite{R16}; there are many
more papers referred to in these works, and numerous applications to
nonlinear problems.

To the question of regularity of solutions,  Ros-Oton and
Serra \cite{RS14a} obtained for (1.5) with $P=(-\Delta )^a$, $0<a<1$, that $f\in L_\infty
$ implies $u\in d^a
C^\alpha (\overline\Omega )$ for small $\alpha $; here
$d(x)=\operatorname{dist}(x,\partial\Omega )$. This was improved later to
$\alpha <a$, and the result was extended to more general
translation-invariant singular integral operators with even kernels
\cite{RS16}.
 
In very recent years, $\psi $do methods have come into the picture
again, mainly through works of the author \cite{G14, G15, G16}, based on the
$a$-transmission property introduced by H\"ormander \cite{H65, H85},
combined with methods from the Boutet de Monvel calculus \cite{B71,
G90, G96, S01, G09}. This approach allows $x$-dependent operators. The
systematic treatment of $\psi $do's primarily takes place in smooth
situations.

It was shown in \cite{G15} that the domain space for solutions of
(1.5) with $f\in L_p(\Omega )$ ($1<p<\infty $, $\Omega $ bounded
smooth), equals $H_p^{a(2a)}(\comega)$, where $H_p^{a(s)}(\comega) =\Lambda
_+^{(-a)}e^+\ol H_p^{s-a}(\Omega )$ (definitions are recalled in Section
2). For a nonhomogeneous Dirichlet problem, one must pass to the larger
space $H_p^{(a-1)(s)}(\comega) $, where the Dirichlet trace operator
$$
\gamma _0^{a-1}u=\Gamma (a)\gamma _0(u/d^{a-1})\tag1.6
$$
(denoted $\gamma _{a-1,0}u$ in \cite{G14, G15, G16}) makes good sense
for $s>a-1/p'$. We generally denote $\gamma _ju=(\partial_n^ju)|_{\partial\Omega }$. 
 Also higher-order trace operators are defined there, e.g.\ the so-called Neumann trace for $s>a+1/p$,
$$
\gamma _1^{a-1}u=\Gamma (a+1)\gamma _1(u/d^{a-1}).\tag1.7
$$
If $\gamma _0^{a-1}u=0$ then $u\in H_p^{a(s)}(\comega)$ and $\gamma
_1^{a-1}u=\gamma _0^au$. Integration-by-parts formulas have up to
now been shown 
%on solution spaces for the homogeneous Dirichlet
%problem (contained in   $H^{a(s)}(\comega)$, $s>a+\frac12$,
%when $\Omega $ is smooth), 
involving
only first boundary values
%$\gamma _0^au$ 
\cite{RS14b, A15, G16} (Corollary 4.5 below); 
they lead
to a useful Pohozaev formula. 

We can now for the first time show a full Green's formula
for $u,v\in H^{(a-1)(s)}(\comega) $ (Theorem 4.4 below): 
%($p=2$):%, $s>a+\frac12$):% $s>a-\max\{1/p,1/p'\}$:
$$
\int _\Omega (Pu\,\bar v-u\,\overline {P^*v})\,dx=\int_{\partial\Omega
}(s_0\gamma _1^{a-1}u\,\gamma _0^{a-1}\bar v-s_0\gamma _0^{a-1}u\, \gamma
_1^{a-1}\bar v
+B\gamma _0^{a-1}u\,\gamma _0^{a-1}\bar v)\, d\sigma 
;\tag1.8
$$
here $s_0(x)$ is the value  of the
principal symbol of $P$ on the interior normal $\nu (x)$ at $x\in\partial\Omega $, $s_0(x)=p_0(x,\nu (x))$, and 
$B$ is a first-order $\psi $do on $\partial\Omega $ (4.20), determined from
the principal and subprincipal symbols of $P$.
To our knowledge, a formula involving both nontrivial Neumann
traces and nontrivial Dirichlet traces has not been found
before. Ellipticity  of $P$ is not assumed.

Next, based on the factorization results of \cite{G16}, we establish
in the elliptic case a detailed formula for the Poisson-like operator
$K_D\colon \varphi \mapsto u$ solving the Dirichlet problem with nonhomogeneous boundary condition
$$
r^+Pu=0\text{ in }\Omega ,\quad \gamma _0^{a-1}u=\varphi \text{ on }\partial\Omega ,\tag1.9
$$
in a parametrix sense, showing exactly how it arises from the
plus-factor $P^+$ in $P\sim P^-P^+$ (Theorems 5.2, 5.4).

As an application, we describe the Dirichlet-to-Neumann operator
(Theorem 6.1),
$$
S_{DN}=\gamma _1^{a-1}K_D;\tag1.10
$$
it is a first-order $\psi $do on $\partial\Omega $ with 
symbol derived from $P^+$. Its principal symbol is, in local coordinates,
$$
s_{DN,0}(x',\xi ')=-\lim_{z_n\to 0+}\Cal F^{-1}_{\xi _n\to
z_n}\log(s_0^{-1}p_0(x',0,\xi )|\xi |^{-2a}) - a|\xi '|.\tag1.11
$$
($S_{DN}$ has only been described before in the elementary case of
$P=(1-\Delta )^a$ on $\rnp$, in \cite{G14}.)
This implies a concrete characterization of the operators $P$ for which the Neumann
problem has Fredholm solvability (Theorem 6.2), namely those for which $S_{DN}$ is
elliptic ($s_{DN,0}$ is nonvanishing for $\xi '\ne 0$). Their
parametrices are likewise described. 

In view of  \cite{G14, G15}, the parametrices and the
Dirichlet-to-Neumann operator we have described 
act in both $H^s_p$-spaces and in
more general Besov and Triebel-Lizorkin spaces; in particular in the
H\"older-Zygmund spaces $C^s_*$ that are of special interest for
nonlinear applications (Corollary 6.4).
\medskip

{\it Plan of the paper:} Section 2 gives preliminaries and notation, and the
Appendix accounts for some results from the Boutet de Monvel calculus
and \cite{G16} that
we use. In Section 3 we show Green's formula for the case of $P=(1-\Delta
)^a$ on $\rnp$, using only Fourier transformation and distribution
theory. In Section 4, Green' formula is shown for general operators $P$
and general smooth domains $\Omega $; this is based on reductions to
applications of the Boutet de Monvel calculus and delicate localization
techniques. Section 5 gives the construction of $K_D$, also using such
tools. In Section 6, $S_{DN}$ is derived, and applied to the discussion
of solvable Neumann problems.

\medskip

\subhead 2. Preliminaries \endsubhead

Multi-index notation is used for
differentiation (and polynomials):
$\partial=(\partial_1,\dots,\partial_n)$, and $\partial^\alpha
=\partial_1^{\alpha _1}\dots \partial_n^{\alpha _n}$ for $\alpha
\in{\Bbb N}_0^n$, with $|\alpha |=\alpha _1+\cdots+\alpha _n$, $\alpha
!=\alpha _1!\dots\alpha _n!$. $D=(D_1,\dots,D_n)$ with $D_j=-i\partial_j$.
Here ${\Bbb N}_0$ denotes the nonnegative ibtegers $\{0,1,2,\dots\}$.
The function $\ang\xi $ stands for $(1+|\xi |^2)^\frac12 $, and $[\xi
]$ denotes a positive
$C^\infty $-function equal to $|\xi |$ for $|\xi |\ge 1$ and $\ge
\tfrac12$ for all $\xi $.

Operators are considered acting on functions or distributions on
${\Bbb R}^n$, and on subsets  
 $\rnpm=\{x\in
{\Bbb R}^n\mid x_n\gtrless 0\}$ (where $(x_1,\dots, x_{n-1})=x'$), and
 bounded $C^\infty $-subsets $\Omega $ with  boundary $\partial\Omega $, and
their complements.
Restriction from $\R^n$ to $\rnpm$ (or from
${\Bbb R}^n$ to $\Omega $ resp.\ $\complement\comega$) is denoted $r^\pm$,
 extension by zero from $\rnpm$ to $\R^n$ (or from $\Omega $ resp.\
 $\complement\comega$ to ${\Bbb R}^n$) is denoted $e^\pm$. Restriction
 from $\crnp$ or $\comega$ to $\partial\rnp$ resp.\ $\partial\Omega $
 is denoted $\gamma _0$. $\Cal S(\crnp)$ stands for $r^+\Cal S({\Bbb
 R}^n)$, where $\Cal S({\Bbb R}^n)$ is Schwartz' space of rapidly
 decreasing $C^\infty $-functions (with dual space $\Cal S'({\Bbb
 R}^n)$, the temperate distributions).

We denote by $d(x)$ a function of the form $
d(x)=\operatorname{dist}(x,\partial\Omega )$ for $x\in\Omega $, $x$ near $\partial\Omega $,
extended to a smooth positive function on $\Omega $; $d(x)=x_n$ in the
case of $\rnp$. Then we define the spaces
$$
\Cal E_\mu (\comega)=e^+\{u(x)=d(x)^\mu v(x)\mid v\in C^\infty
(\comega)\},\tag2.1
$$
for $\operatorname{Re}\mu
 >-1$; for other $\mu $, cf.\ \cite{G15}.

A {\it pseudodifferential operator} ($\psi $do) $P$ on ${\Bbb R}^n$ is
defined from a symbol $p(x,\xi )$ on ${\Bbb R}^n\times{\Bbb R}^n$ by 
$$
Pu=p(x,D)u=\operatorname{OP}(p(x,\xi ))u 
=(2\pi )^{-n}\int e^{ix\cdot\xi
}p(x,\xi )\hat u\, d\xi =\Cal F^{-1}_{\xi \to x}(p(x,\xi )\hat u(\xi
));\tag 2.2
$$  
here $\Cal F$ is the Fourier transform $(\F u)(\xi )=\hat u(\xi
)=\int_{{\Bbb R}^n}e^{-ix\cdot \xi }u(x)\, dx$. 
We refer to
textbooks such as H\"ormander \cite{H85}, Taylor \cite{T91}, Grubb \cite{G09} for the rules of
calculus. \cite{G09} moreover gives an account of the Boutet de Monvel
calculus of {\it pseudodifferential boundary problems}, cf.\ also e.g.\
\cite{G96, S01}. 

We take $p$ in the symbol space $S^m_{1,0}({\Bbb R}^n\times{\Bbb R}^n)$, consisting of
$C^\infty $-functions $p(x,\xi )$
such that $\partial_x^\beta \partial_\xi ^\alpha p(x,\xi
)$ is $O(\ang\xi ^{m-|\alpha |})$ for all $\alpha ,\beta $, for some
$m\in{\Bbb R}$ (global estimates); then $p$ and $P$ have order $m$.
A symbol $p$ is said to be {\it classical} when it moreover 
has an asymptotic expansion $p(x,\xi )\sim \sum_{j\in{\Bbb
N}_0}p_j(x,\xi )$ with $p_j$ homogeneous in $\xi $ of degree $m-j$ for
$|\xi |\ge 1$, all $j$, and $p(x,\xi )- \sum_{j<J}p_j(x,\xi )\in S^{m-J}_{1,0}({\Bbb
R}^n\times \R^n)$ for all $J\in {\Bbb N}_0$.

Recall in particular the composition rule: When $PQ=R$, then $R$ has
a symbol $r(x,\xi )$ with the following asymptotic expansion, called the
Leibniz product:
$$
r(x,\xi )\sim p(x,\xi )\# q(x,\xi )= {\sum}_{\alpha \in{\Bbb
N}_0^n}\partial_\xi ^\alpha p(x,\xi ) D_x^aq(x,\xi )/\alpha !.\tag 2.3
$$

When $P$ is a $\psi $do on ${\Bbb R}^n$, $P_+=r^+Pe^+$
denotes its truncation to $\rnp$, or to $\Omega $, depending on the context.
%\medskip

Let $1<p<\infty $  (with $1/p'=1-1/p$), then the $L_p$-Sobolev spaces
(Bessel-potential spaces) are defined for $s\in{\Bbb R}$ by
$$
\aligned
H^s_p(\R^n)&=\{u\in \SD'({\Bbb R}^n)\mid \F^{-1}(\ang{\xi }^s\hat u)\in
L_p(\R^n)\},\\
\dot H^{s}_p(\comega)&=\{u\in H^{s}_p({\Bbb R}^n)\mid \supp u\subset
\comega \},\\
\ol H^{s}_p(\Omega)&=\{u\in \D'(\Omega )\mid u=r^+U \text{ for some }U\in
H^{s}_p(\R^n)\};
\endaligned 
$$
here $\operatorname{supp}u$ denotes the support of $u$. The definition
is also used with $\Omega =\rnp$. In most current texts, $\ol
H^s_p(\Omega )$ is denoted $H^s_p(\Omega )$ without the overline (that
was introduced along with the notation $\dot H$ in \cite{H65, H85}), but we keep it here since it is
practical in indications of dualities, and makes the notation more
clear in formulas where
both types occur. When $p=2$, the mention of
$p$ is left out.
We recall that $\ol H_p^s(\Omega )$ and $\dot H_{p'}^{-s}(\comega)$ are dual
spaces with respect to a sesquilinear duality extending the $L_2(\Omega )$-scalar
product, written e.g.\
$$
\ang{f,g}_{\ol H_p^s(\Omega ),\dot H_{p'}^{-s}(\comega)},\text{ or just }\ang{f,g}_{\ol H_p^s,\dot H_{p'}^{-s}}.
$$

There are many other interesting scales of spaces,
the
Triebel-Lizorkin spaces $F^s_{p,q}$ and   Besov spaces  $B^s_{p,q}$
($B^s_p=B^s_{p,p}$), where 
the problems can be studied; see details in \cite{G14}. This includes the H\"older-Zygmund spaces $B^s_{\infty
,\infty }$, also denoted $C^s_*$; they are interesting because $C^s_*({\Bbb R}^n)$ equals
the H\"older space $C^s({\Bbb R}^n)$ when $s\in\rp\setminus {\Bbb
N}$. 
There are also local versions: $H^s_{p,\operatorname{loc}}(\Omega )$
consists of the distributions $u\in \Cal D'(\Omega )$ such that $\varphi u\in \ol
H^s_p({\Bbb R}^n)$ for all $\varphi \in C_0^\infty (\Omega )$, and $
H^s_{p,\operatorname{comp}}(\Omega )$ consists of the elements of $
H^s_{p}({\Bbb R}^n)$ with compact support in $\Omega $. There are similar definitions where $H^s_p$ is replaced by
$F^s_{p,q}$ or  $B^s_{p,q}$.

A classical $\psi $do $P$ of order $m\in{\Bbb R}$ maps $
H^s_{p,\operatorname{comp}}({\Bbb R}^n )$ into $
H^{s-m}_{p,\operatorname{loc}}({\Bbb R}^n )$ for all $s\in{\Bbb R}$. It is {\it elliptic}, when
$p_0(x,\xi )\ne 0$ for $|\xi |\ge 1$; then for any $u\in \Cal E'({\Bbb
R}^n)$, any open $\Omega $,
$Pu\in H^{s-m}_{p,\operatorname{loc}}(\Omega )$ implies $u\in
H^{s}_{p,\operatorname{loc}}(\Omega )$.
 Analogous
results hold in the other scales  $F^s_{p,q}$ and $B^s_{p,q}$.

H\"ormander introduced in
\cite{H65, H85} {\it the $\mu
$-transmission condition} at $\partial\Omega $ for a classical $\psi
$do of order $m$, $\mu
\in{\Bbb C}$: In local coordinates,
$$
\partial_x^\beta \partial_\xi ^\alpha p_j(x,-\nu )=e^{\pi i(m-2\mu -j-|\alpha | )
}\partial_x^\beta \partial_\xi ^\alpha p_j(x,\nu ),\tag 2.4
$$
for all $x\in\partial\Omega $, all $j,\alpha ,\beta $, where 
$\nu $ denotes the interior normal to $\partial\Omega $ at $x$.
The Boutet de Monvel calculus treats the case where $\mu =0$, and the
operators $P$ we shall study satisfy it with $\mu =a>0$.

A special role in the theory is played by the {\it order-reducing
operators}. There is a simple definition of operators $\Xi _\pm^t $ on
${\Bbb R}^n$, $t\in{\Bbb R}$,
$$ 
\Xi _\pm^t =\operatorname{OP}(\chi _\pm^t),\quad \chi _\pm^t=(\ang{\xi '}\pm i\xi _n)^t ;\tag 2.5 
$$
 they preserve support
in $\crnpm$, respectively. The functions
$(\ang{\xi '}\pm i\xi _n)^t $ do not satisfy all the estimates
required for the class $S^{t }_{1,0}({\Bbb
R}^n\times{\Bbb R}^n)$, but the operators are useful for many
purposes. There is a more refined choice $\Lambda _\pm^t $
\cite{G90, G15}, with
symbols $\lambda _\pm^t (\xi )$ that do
satisfy all the estimates for $S^{ t }_{1,0}({\Bbb
R}^n\times{\Bbb R}^n)$; here $\overline{\lambda _+^t }=\lambda _-^{t }$.
The symbols have holomorphic extensions in $\xi _n$ to the complex
halfspaces ${\Bbb C}_{\mp}=\{z\in{\Bbb C}\mid
\operatorname{Im}z\lessgtr 0\}$, and hence the operators preserve
support in $\crnpm$, respectively; operators with that property are
called "plus" resp.\ "minus" operators. There is also a pseudodifferential definition $\Lambda
_\pm^{(t )}$ adapted to the situation of a smooth domain $\Omega
$, cf.\ \cite{G15}.

It is elementary to see by the definition of the spaces $H^s_p(\R^n)$
in terms of Fourier transformation, that the operators define homeomorphisms 
for all $s$:
$$
\Xi^t _\pm\colon H^s_p(\R^n) \simto H^{s- t
}_p(\R^n), \quad  
\Lambda ^t _\pm\colon H^s_p(\R^n) \simto H^{s- t
}_p(\R^n)\tag 2.6
$$
(and so does of course $\Xi ^t =\Op(\ang \xi ^t )=\ang D^t$). The special
interest is that the "plus"/"minus" operators also 
 define
homeomorphisms related to $\crnp$ and $\comega$, for all $s\in{\Bbb R}$: 
$$
\aligned
\Xi ^{t }_+\colon \dot H^s_p(\crnp )\simto
\dot H^{s- t }_p(\crnp),&\quad
r^+\Xi ^{t }_{-}e^+\colon \ol H^s_p(\rnp )\simto
\ol H^{s- t }_p(\rnp );
\\
\Lambda^{(t )}_+\colon \dot H^s_p(\comega )\simto
\dot H^{s- t }_p(\comega ),&\quad
r^+\Lambda ^{(t )}_{-}e^+\colon \ol H^s_p(\Omega )\simto
\ol H^{s- t }_p(\Omega );
\endaligned \tag 2.7
$$
 here %$\Xi ^t _{-,+}$ resp.\  $\Lambda ^{(t )}_{-,+}$ is short for
$r^+\Xi ^t _-e^+$  
resp.\  $r^+\Lambda ^{(t )}_{-}e^+$ are suitably extended to large negative $s$ (cf.\ Rem.\ 1.1
 and Th.\ 1.3 in \cite{G15}). 
The first line in (2.7) also holds with
 $\Xi $ replaced by $\Lambda $.

One has moreover that
 the operators $\Xi ^t _{+}$ and $r^+\Xi ^{t }_{-}e^+$ identify with each other's adjoints
over $\crnp$, because of the support preserving properties; more precisely,
$$%\aligned
\Xi ^{t }_{+ }\colon \dot H_{p'}^{ t -s}(\crnp)\to \dot
 H_{p'}^{-s}(\crnp)  \text{ and }r^+\Xi ^{t }_{-}e^+\colon \ol H_p^s(\rnp)\to
\ol H_p^{s- t }(\rnp)\text{ are adjoints},
%\endaligned
\tag 2.8
$$
for  all $s\in{\Bbb R}$. The same holds
for the operators $\Lambda _+^t ,r^+\Lambda ^{t }_{-}e^+$, and there is a
similar statement for $\Lambda ^{(t )}_+$ and $r^+\Lambda ^{(
t )}_{-}e^+$ relative to the set $\Omega $.

Now we rapidly recall the features of some special spaces studied
in \cite{G15}; a more detailed account of some particular instances is
given below in Section 3. They are the
{\it $\mu $-transmission spaces} 
introduced by
H\"ormander \cite{H65} (for $p=2$), cf.\ \cite{G15}, which are particularly adapted to $\mu
$-transmission operators $P$ (we just take real $\mu >-1$): 
$$
\aligned
H^{\mu (s)}_p(\crnp)&=\Xi _+^{-\mu }e^+\ol H_p^{s- \mu
}(\rnp)=\Lambda  _+^{-\mu }e^+\ol H_p^{s- \mu
}(\rnp)
,\quad  s> \mu -1/p',\\
H^{\mu (s)}_p(\comega)&=\Lambda  _+^{(-\mu )}e^+\ol H_p^{s- \mu
}(\Omega ),\quad  s> \mu -1/p'.
\endaligned\tag 2.9
$$
In fact, $r^+P$ (of order $m$) maps them into %$C^\infty (\comega)$,  
$\ol H_p^{s- m}(\rnp)$ resp.\ $\ol
 H_p^{s- m}(\Omega)$ 
(cf.\
 \cite{G15} Sections 1.3, 2, 4),
 and they represent, when $P$ is
 elliptic, the solution space for the homogeneous Dirichlet problem
 (1.5) with $f\in \ol H^{s-m}_p(\rnp)$ resp.\ $\ol H^{s-m}_p(\Omega)$.
Moreover, $r^+P$  maps $\Cal E_\mu
 (\comega)$ into 
$C^\infty (\comega)$, and $\Cal E_\mu (\comega)$ is the solution space for the
Dirichlet problem with data in $C^\infty (\comega)$.  $\Cal E_\mu (\comega)$ is dense in  $H^{\mu
(s)}_p(\comega)$ for all $s$, and $\bigcap_s H^{\mu
(s)}_p(\comega)=\Cal E_\mu (\comega)$. (For $\Omega =\rnp$, $\Cal
E_\mu (\crnp)\cap \Cal E'({\Bbb R}^n)$ is dense in $H^{\mu
(s)}_p(\crnp)$ for all $s$.)

 The following trace operators $\gamma _k^\mu $ (denoted $\gamma _{\mu
,k}$ in \cite{G14, G15, G16}), are defined on $\Cal E_\mu (\comega)$,
and on $H_p^{\mu (s)}(\comega)$ for $s>\mu +k+1/p$:
$$
\gamma _k^\mu u=c_k\gamma _k(u/d^{\mu }), \quad c_k=\Gamma (\mu +1+k),
\tag 2.10
$$
where the $\gamma _k$ are the standard trace operators  $\gamma
_ku=(\partial_n^ku)|_{\partial\Omega }$, $k=0,1,\dots$. Here when
$M\in{\Bbb N}$ and $s>\mu +M-1/p'$, 
$\varrho ^\mu _M=\{\gamma ^\mu _0, \dots,\gamma ^\mu
_{M-1}\}$ maps surjectively
$$
\varrho ^\mu _M\colon H^{\mu (s)}_p(\comega)\to {\prod}_{j=0}^{M-1}B^{s-\mu
-j-1/p}_p(\partial\Omega ),\text{ with kernel } H^{(\mu +M-1)(s)}_p(\comega),
$$
cf.\ \cite{G15} Th.\ 5.1.

 One has that $H^{\mu
 (s)}_p(\comega)\supset \dot H_p^s(\comega)$, and that the distributions are
 locally in $H^s_p$   on $\Omega $, but at the boundary they in general have a
 singular behavior (cf.\ \cite{G15} Th.\ 5.4):
$$
H_p^{\mu (s)}(\comega)\cases =\dot H_p^s(\comega) \text{ if }s\in
\,]\mu -1/p',\mu +1/p[\,,\\
\subset e^+d^\mu  \ol H_p^{s- \mu }(\Omega)+\dot
H_p^{s}(\comega)\text{ if }s>\mu +1/p, \; s-\mu -1/p\notin {\Bbb N}.
\endcases
\tag 2.11
$$
The inclusion in the second line is not an identity, but elements of
$e^+d^\mu  \ol H_p^{s- \mu }(\Omega)$ enter nontrivially. In local
coordinates where $\Omega $ is replaced by $\rnp$, there are elements
of the form $e^+x_n^\mu K_0\varphi $, where $K_0$ is the standard
Poisson operator $K_0\varphi =\F^{-1}_{\xi '\to x'}[r^+e^{-\ang{\xi
'}x_n}\hat\varphi (\xi ')]$, and $\varphi $ runs through $B^{s-\mu
-1/p}_p({\Bbb R}^{n-1})$. In particular, there are elements lying
in   $e^+d^\mu \ol H_p^{s- \mu
}(\Omega)$, not in $e^+d^\mu \ol H_p^{s'- \mu }(\Omega)$ for any  $s'>s$.

When $-1<\mu <0$, the factor $d^{\mu }$ blows up at the boundary, so
functions in these spaces can be viewed as ``large''. However, it is
worth remarking that they are  $L_p$-integrable for $p>-1/\mu $, $s\ge 0$. For
$0\le s<\mu +1/p$ this follows since $H_p^{(\mu ) (s)}(\comega)=\dot
H_p^s(\comega)$ then, and for larger $s$,
it follows since $H_p^{(\mu ) (s')}(\comega)\subset H_p^{(\mu )
(s)}(\comega)$ when $s'>s$.

In the present paper, we use these spaces with $\mu =a$ and with $\mu
=a-1$, mainly with $p=2$. For example, 
the spaces $H^{(a-1)(s)}(\comega)$ satisfy
$$
H^{(a-1) (s)}(\comega)\cases =\dot H^s(\comega) \text{ if }s\in
\,]a-\frac32,a-\frac12[\,,\\
\subset e^+d^{a-1}  \ol H^{s- a+1 }(\Omega)+\dot
H^{s}(\comega)\text{ if }s>a-\frac12,\; s-a-\frac12 \notin {\Bbb N}.
\endcases
\tag 2.12
$$
The first two traces $\gamma _0^{a-1}u$ (Dirichlet) and $\gamma
_1^{a-1}u$ (Neumann) are well-defined when $s>a-\frac12$ resp.\
$s>a+\frac12$. Here  
$u\in H^{(a-1) (s)}(\comega)$ with $\gamma _0^{a-1} u=0$ implies $u\in
H^{a (s)}(\comega)$.

Our Green's formula will be shown for $p=2$, since it is closely
connected with $L_2$-scalar products. The calculations of solution
operators will be performed with $p=2$, and
supplied with remarks on how the results extend to $H^s_p$- and other spaces as in
\cite{G14, G15}. 
%$p\ne 2$, using
%that $\Cal E_\mu (\comega)$ is dense in $H_p^{\mu  (s)}(\comega)$. 

Further prerequisites are collected in the Appendix.

\subhead 3. Green's formula in the simplest case \endsubhead

We begin by an elementary explanation of the Dirichlet and Neumann
boundary values on $\rnp$. 
Consider $u,v\in \Cal E_{a-1}(\crnp)\cap e^+\Cal S(\crnp)$, with
$a>0$. %(Here $\Cal S(\crnp)$ is short for $r^+\Cal S({\Bbb R}^n)$,
%also denoted $\Cal S_+$. 
The
boundary values $\gamma _k^{a-1}u=u_k$ are defined from the expansion
$$
u(x)=u_0(x')I^{a-1}(x_n)+u_1(x')I^{a}(x_n)+\dots+ u_k(x')I^{a-1+k}(x_n)+O(x_n^{a+k}),\tag3.1
$$
where (as in \cite{G15}) $I^\mu (x_n)=H(x_n)x_n^\mu /\Gamma (\mu +1)$
when $\operatorname{Re}\mu >-1$, $H$ is the Heaviside function. The
Gamma factor serves to normalize $I^\mu
$  so that $\partial_{x_n}I^\mu =I^{\mu -1}$; this formula is also
used to define the distribution for lower $\operatorname{Re}\mu $. 
The expansion (3.1) follows from a Taylor expansion of $w(x)=u(x)/x_n^{a-1}$ for
$x_n\to 0+$:
$$
u(x)=x_n^{a-1}w(x',0)+x_n^{a}\partial_nw(x',0)+\dots +x_n^{a+k-1}\tfrac1{k!}\partial_n^kw(x',0)+O(x_n^{a+k}).\tag3.2
$$
Note in particular that 
$$
\aligned
\gamma _0^{a-1}u&=u_0=\Gamma (a)\gamma _0w=\Gamma (a)\gamma
_0(u(x)/x_n^{a-1}),\\
 \gamma _1^{a-1}u&=u_1=\Gamma (a+1)\gamma _1w=\Gamma (a+1)\gamma
 _1(u(x)/x_n^{a-1});
\endaligned
\tag3.3
$$
they will be
viewed as the {\it Dirichlet 
resp.\ Neumann
traces} of $u\in \Cal E_{a-1}(\crnp)$. 

\example{Remark 3.1}
When $a=1$, i.e., $u\in e^+C^\infty (\crnp)$, this fits together with the
usual convention for Dirichlet and Neumann traces associated with the
Laplacian. Note however that for large $a$, the names are  used here
for the two lowest nontrivial traces. E.g., if $u\in \E_k(\crnp)$ for
a large integer
$k$, whereby
$u=x_n^{k}v$ for a $v\in C^\infty (\crnp )$, then the first $k$ standard traces of $u$ vanish, and $\{\gamma ^k_0u,\gamma ^k_1u\}=\{k!\gamma
_0v,(k+1)!\gamma _1v\}$. For $2k$-order elliptic differential
operators there is another convention for Dirichlet and Neumann values as
the trace collections $\{\gamma _0u,\dots,\gamma _{k-1}u\}$ resp.\
$\{\gamma _ku,\dots,\gamma _{2k-1}u\}$. Also for fractional operators,
there are well-posed boundary value problems with sets of traces, as
e.g.\ in \cite{G15} Th.\ 6.1.
\endexample

Besides the expansion (3.1) it will be convenient to use an
expansion where the partially Fourier transformed 
terms have a factor $e^{-\sigma x_n}$, $\sigma =\ang{\xi '}$:
$$
\aligned
\Cal F_{x'\to\xi '}u&=\acute u(\xi ',x_n)
=\hat \varphi _0(\xi ')I^{a-1}(x_n) e^{-\sigma x_n}+\acute u'(\xi ',x_n)\\
&=\hat
\varphi _0(\xi ')I^{a-1}(x_n) e^{-\sigma x_n}+\hat \varphi _1(\xi ')I^{a}(x_n) e^{-\sigma
x_n}+\acute u''(\xi ',x_n)\\
&\sim\sum_{k\ge 0}\hat \varphi _k(\xi ')I^{a-1+k}(x_n)e^{-\sigma
x_n};\text{ here}\\
\Cal Fu&\sim \sum_{k\ge 0}\hat \varphi _k(\xi ')(\sigma +i\xi _n)^{-a-k},
\endaligned\tag3.4
$$
using the formula 
$$
\Cal F_{x_n\to \xi _n}[I^\mu(x_n) e^{-\sigma x_n}]=(\sigma +i\xi _n)^{-\mu
-1}.\tag3.5
$$

Correspondingly,
$$
%\aligned
u=U_0+u'=U_0+U_1+u''\sim \sum_{k\ge 0}U_k,\text{ with }U_k=\Cal F^{-1}_{\xi '\to x'}[\hat \varphi _k(\xi ')I^{a-1+k}(x_n)e^{-\sigma
x_n}].%\
 %&=V_0+v'=V_0+V_1+v'',\text{ with }V_k=\Cal F^{-1}_{\xi '\to x'}[ v_k(\xi ')I^{a-1+k}e^{-\sigma
%_n}].
%endaligned
\tag3.6
$$
Here we note that $u'\in \Cal E_a$ and $u''\in \Cal E_{a+1}$, with
$U_0\in \Cal E_{a-1}$, $U_1\in \Cal E_a$. So $\gamma _0^{a-1}u'=0$
with $\gamma _0^{a-1}u=\gamma _0^{a-1}U_0$, and $\gamma _0^{a}u''=0$
with $\gamma _0^{a}u'=\gamma _0^{a}U_1$. (They are all in $e^+\Cal S(\crnp)$.) 
%Since the $\varphi _k$ are in $\Cal S({\Bbb R}^{n-1})$, $u'\in \E_{a}(\crnp)\ca%p e^+\Cal S(\crnp)$ and $u''\in \E_{a+1}(\crnp)\cap
%e^+\Cal S(\crnp)$. 

The transition between the coefficient sets $\{u _k\}$ and  $\{\varphi
_k\}$ can be found by comparison of
the expansion of  $\Cal F_{x'\to\xi '}w=\acute w(\xi ',x_n)$ %(cf.\ (3.2a)) 
with the
expansion  we get from $\acute w_e(\xi
',x_n)=e^{\sigma x_n}\acute w(\xi
',x_n)=e^{\sigma x_n}\acute u(\xi ',x_n)/x_n^{a-1}$.
Let us just do this in detail for the first two coefficients
(sufficient for the present paper):
$$
\aligned
e^{\sigma x_n}\acute u&=x_n^{a-1}\acute w_e(\xi
',0)+x_n^a\partial_n\acute w_e(\xi ',0)+O(x_n^{a+1})\\
&=e^{\sigma x_n}[x_n^{a-1}\hat
w_0(\xi ')+x_n^a(\hat w_1(\xi ')+\sigma \hat w_0(\xi
'))+O(x_n^{a+1})].
\endaligned
$$
We see that $\hat\varphi _0=\Gamma (a)\hat w_0=\hat u_0$ and $\hat \varphi _1=\Gamma
(a+1)(\hat w_1+\sigma \hat w_0)=\hat u_1+a\sigma \hat u_0$, so
$$
\varphi _0=u_0,\quad \varphi _1=u_1+a\ang{D'}u_0,\tag3.7
$$
where $\ang{D'}=\OP(\ang{\xi '})$.
This allows us to relate the functions $u_0,u_1$ to boundary values of $\Xi
_+^{a-1}u$, cf.\ (2.5). In view of (3.4), $\Xi _+^{a-1}=\operatorname{OP}((\sigma
+i\xi _n)^{a-1})$ has the effect
$$
\Xi _+^{a-1}u\sim \Cal F^{-1}\sum_{k\ge 0}\hat \varphi _k(\xi ')(\sigma +i\xi _n)^{-1-k}.
$$
Then since $\gamma
_0I^k=0$ for $k=1,2,\dots$ (recall that the trace $\gamma _0$  is taken from $\rnp$), 
$$
\aligned
\gamma _0\Xi _+^{a-1}u&=\gamma _0\Cal F^{-1}(\hat \varphi _0(\xi ')(\sigma +i\xi
_n)^{-1}+\hat \varphi _1(\xi ')(\sigma +i\xi _n)^{-2}+\dots)\\
&=\gamma _0\Cal F^{-1}_{\xi '\to x'}(\hat \varphi _0(\xi ')I^0e^{-\sigma
x_n}+\hat \varphi _1(\xi ')I^1e^{-\sigma x_n}+\dots)\\
&=\Cal F^{-1}_{\xi '\to
x'}\hat \varphi _0(\xi ')=\varphi _0=u_0,\\
\gamma _0\partial_n\Xi _+^{a-1}u
%&=\gamma _0\partial_n\Cal F^{-1}(\hat \varphi _0(\xi ')(\sigma +i\xi
%_n)^{-1}+\hat \varphi _1(\xi ')(\sigma +i\xi _n)^{-2}+\dots)\\
&=\gamma _0\partial_n\Cal F^{-1}_{\xi '\to x'}(\hat \varphi _0(\xi ')I^0e^{-\sigma
x_n}+\hat \varphi _1(\xi ')I^1e^{-\sigma x_n}+\dots)\\
&=\Cal F^{-1}_{\xi '\to x'}(-\hat \varphi _0(\xi ')\sigma +\hat \varphi _1(\xi '))
=-\ang{D'}\varphi _0 +\varphi _1=u_1+(a-1)\ang{D'}u_0,
\endaligned
$$
cf.\ (3.7), showing that 
$$
u_0=\gamma _0\Xi _+^{a-1}u,\quad u_1=\gamma _0\partial_n\Xi _+^{a-1}u-(a-1)\ang{D'}u_0.\tag3.8
$$
(It is used that $\Cal F^{-1}_{\xi '\to x'}(\hat \varphi _0(\xi ')I^{-1}e^{-\sigma
x_n})= \varphi _0(x ')\delta (x_n)$ does not
contribute to the boundary value from $\rnp$.) (3.8) was
also shown in \cite{G15} Sect.\ 5;
related calculations occur in \cite{G14}, Appendix. 

For $u'=u-U_0$ (cf.\ (3.6)), we note that since $u'\in \E_a$ with the expansion $U_1+u''$, $u''\in \E_{a+1}$, 
$$
\gamma _0^au'=\varphi _1=\gamma _1^{a-1}u+a\ang{D'}\gamma _0^{a-1}u=u_1+a\ang{D'}u_0.\tag3.9
$$
In particular, $u'$ itself satisfies, since $\gamma _0^{a-1}u'=0$ and
$\gamma _1^{a-1}u'=\gamma _1^{a-1}u$,
$$
\gamma _0^au'=\gamma _1^{a-1}u'.\tag3.10
$$
In other words: {\it When $u\in  \Cal E_{a-1}$ %(more generally $u\in
%H^{(a-1)(s)}(\crnp)$) 
is such that the Dirichlet trace $\gamma
_0^{a-1}u$ of $u$
vanishes, then the
Neumann trace equals $\gamma _0^au$.}

Let $v$  be another function in $ \Cal E_{a-1}(\crnp)\cap e^+\Cal
S(\crnp)$; then we expand it similarly as in (3.6) with
coefficients $\hat \psi _k$:$$
%\aligned
v=V_0+v'=V_0+V_1+v'',\text{ with }V_k=\Cal F^{-1}_{\xi '\to x'}[\hat \psi _k(\xi ')I^{a-1+k}e^{-\sigma
x_n}].%\
 %&=V_0+v'=V_0+V_1+v'',\text{ with }V_k=\Cal F^{-1}_{\xi '\to x'}[ v_k(\xi ')I^{a-1+k}e^{-\sigma
%_n}].
%endaligned
\tag3.11
$$
%The statement after (3.10) extends to these spaces.

The formula (3.8) allows us to deduce mapping properties of the
$\gamma _k^{a-1}$ in Sobolev spaces. Recall from Section 2 that $\Cal
E_{a-1}(\crnp)\cap \Cal E'({\Bbb R}^n)$ is
dense in $ H^{(a-1)(s)}(\crnp)=\Xi
_+^{1-a}e^+\ol H^{s-a+1}(\rnp)$. Then since $\gamma _k\colon \ol
H^{s-a+1}(\rnp)\to H^{s-a-k+\frac12}({\Bbb R}^{n-1})$, the Dirichlet and
Neumann traces extend by continuity to
continuous operators:
$$
\aligned
\gamma _0^{a-1}&\colon H^{(a-1)(s)}(\crnp)\to
H^{s-a+\frac12}({\Bbb R}^{n-1}),\quad s>a-\tfrac12,\\
\gamma _1^{a-1}&\colon H^{(a-1)(s)}(\crnp)\to
H^{s-a-\frac12}({\Bbb R}^{n-1}),\quad s>a+\tfrac12.
\endaligned\tag3.12
$$
In this context, we note that (since $\varphi _0=u_0$)
$$\aligned
U_0&=\Cal F^{-1}_{\xi '\to x'}[\hat u_0(\xi ')I^{a-1}(x_n)e^{-\sigma
x_n}]
=\Cal F^{-1}_{\xi \to x}[\hat u_0(\xi ')(\sigma +i\xi
_n)^{-a}]\\
&=\Xi _+^{1-a}\Cal F^{-1}_{\xi \to x}[\hat u_0(\xi ')(\sigma +i\xi _n)^{-1}]=\Xi _+^{1-a}e^+K_0u_0,
\endaligned$$
where $K_0$ is the well-known Poisson operator
$
K_0\varphi =\Cal F^{-1}_{\xi '\to x'}[\hat \varphi (\xi ')r^+e^{-\ang{\xi '}
x_n}]
$,
defining a right-inverse of $\gamma _0$; its symbol is $(\sigma +i\xi _n)^{-1}=\chi _+^{-1}$. It maps $H^t({\Bbb
R}^{n-1})\to \ol H^{t+\frac12}(\rnp)$ for all $t\in{\Bbb R}$, so 
$$
K_0^{a-1}\equiv \Xi _+^{1-a}e^+K_0\colon 
H^{s-a+\frac12}({\Bbb R}^{n-1})\to  H^{(a-1)(s)}(\crnp),\quad s>a-\tfrac12,
$$
defining a right-inverse of $\gamma _0^{a-1}$. (This was also shown in
\cite{G15} Cor.\ 5.3.)

The various identities shown above, and the remark after (3.10), extend to these spaces, for $s$ suitably chosen.

We shall now give a relatively elementary proof of the desired Green's formula for
$P=(1-\Delta )^a$ on $\rnp$. Since $(1+|\xi |^2)^a=(\ang{\xi '}-i\xi
_n)^a(\ang{\xi '}+i\xi _n)^a$, $(1-\Delta )^a=\Xi ^a _-\Xi
^a_+$. In view of (3.5):
$$
\aligned
\Xi _+^aU_k&=\Cal F^{-1}[(\sigma +i\xi _n)^a\hat \varphi _k(\xi ')(\sigma +i\xi
_n)^{-a-k}]\\
&=\Cal F^{-1}_{\xi '\to x'}[\hat \varphi _k(\xi
')I^{k-1}e^{-\sigma x_n}]\text{ for }k\in{\Bbb N_0},          \text{
in particular,}\\
\Xi _+^aU_0&=\Cal F^{-1}[(\sigma +i\xi _n)^a\hat \varphi _0(\xi ')(\sigma +i\xi
_n)^{-a}]=\varphi _0(x')\otimes \delta (x_n),
\\
\Xi _+^aU_1&=\Cal F^{-1}[(\sigma +i\xi _n)^a\hat \varphi _1(\xi ')(\sigma +i\xi
_n)^{-a-1}]=\Cal F^{-1}_{\xi '\to x'}[\hat \varphi _1(\xi ')He^{-\sigma x_n}].
\endaligned\tag 3.13
$$
An application of $\Xi
_-^a$ gives
$$
\aligned
(1-\Delta )^aU_0&=\Xi _-^a(\varphi _0(x')\otimes \delta (x_n)),\text{
supported in }\crnm,\\
(1-\Delta )^aU_1&=\Cal F^{-1}[(\sigma -i\xi _n)^a\hat \varphi _1(\xi ')(\sigma +i\xi
_n)^{-1}].
\endaligned
$$
From the first line we conclude for $P=(1-\Delta )^a$:
$$
r^+PU_0=0,\text{ hence }r^+Pu=r^+Pu'.\tag3.14
$$

Now  $\int_{\rnp}r^+Pu\,\bar v\, dx$ will be worked out. 
We know from \cite{G15} Th.\ 4.2 that $r^+P$ maps
$H^{(a-1)(s)}(\crnp)$ into $\ol H^{s-2a}(\rnp)$ (and  $\Cal E_{a-1}(\crnp)\cap \Cal
E'({\Bbb R}^n)$ into $\bigcap_t \overline H^{t}(\rnp)$). For large
$s$, $\ol H^{s-2a}(\rnp)$ is
a space of continuous functions, and  $H^{(a-1)(s)}(\crnp)$ is as such,
supplied with
continuous functions multiplied by  $x_n^{a-1}$, so $r^+Pu\,\bar v$ is
integrable for $x_n\to 0$. For smaller $s$, we need an interpretation
as a duality.

Note first (cf.\ (3.11) and (3.14)) that
$$
\int_{\rnp}r^+Pu\,\bar v\, dx=\int_{\rnp}r^+Pu'\,\bar v\, dx=\int_{\rnp}r^+Pu'\,\bar v'\, dx+\int_{\rnp}r^+Pu'\,\bar V_0\, dx.\tag3.15
$$
In the term $\int_{\rnp}r^+Pu'\,\bar v'\, dx$, $v'\in x_n^{a}\ol
H^{s-a}(\rnp)+\dot H^{s}(\crnp)$ does not give integrability problems. This integral will be
left unchanged, to match a similar integral with $P$ applied to $v'$. It is the last integral that will be
reduced to an integral of boundary values, and which we now study more closely.

From now on,
%and we now pass to these larger spaces, 
take $u$ and $v$ in
$H^{(a-1)(s)}(\crnp)=\Xi
_+^{1-a}e^+ \ol H^{s-a+1}(\rnp)$. 
The
following calculations are very similar to those in the proof of
\cite{G16}, Th.\ 3.1. 
Let
$s>a+\frac12 $, then (for small $\varepsilon >0$)
$$
\aligned
&u,v,U_0,V_0\in H^{(a-1)(s)}(\crnp)= \Xi _+^{1-a}e^+\overline
H^{s-a+1}(\rnp)\subset\Xi _+^{1-a}\dot H^{\frac12-\varepsilon
}(\crnp)=\dot H^{a-\frac12-\varepsilon }(\crnp),\\
&u',v'\in H^{a(s)}(\crnp)= \Xi _+^{-a}e^+\overline
H^{s-a}(\rnp)\subset\Xi _+^{-a}\dot H^{\frac12-\varepsilon
}(\crnp)=\dot H^{a+\frac12-\varepsilon }(\crnp),\\
&r^+Pu, r^+Pu', r^+Pv, r^+Pv' \in \ol H^{s-2a}(\rnp)\subset \ol H^{\frac12-a+\varepsilon }(\rnp),\\
&u_k, v_k\in 
 H^{s-a-k+\frac12}({\Bbb R}^{n-1})\subset H^{1-k+\varepsilon }({\Bbb
 R}^{n-1}),\quad k=0,1.
\endaligned\tag3.16
$$ 
We then find %recalling that $r^+Pu=r^+Pu'$, 
that $I=
\int_{\rnp}r^+Pu'\,\bar V_0\,dx$ can be interpreted in these larger
spaces as
$$
I\equiv\ang{r^+Pu',V_0}_{ \ol H^{\frac12-a+\varepsilon }(\rnp) ,\dot H^{a-\frac12-\varepsilon
}(\crnp)}\tag3.17
$$
(sesquilinear duality). Since $u'\in H^{a(s)}(\crnp)$, $r^+Pu'=r^+\Xi _-^ae^+r^+\Xi
_+^au'=r^+\Xi _-^ae^+ w$,
where $w=r^+\Xi _+^au'\in \ol H^{s-a}(\rnp)$; here $r^+\Xi _-^ae^+$
maps $\ol H^t(\rnp)$ to $\ol H^{t-a}(\rnp)$ for all $t\in{\Bbb R}$,
with adjoint 
$\Xi _+^a\colon \dot H^{a-t}(\crnp)\to  \dot H^{-t}(\crnp)$, cf.\
(2.7), (2.8). 
Therefore, by (3.13),
$$
\aligned
I&=\ang{r^+\Xi _-^ae^+r^+\Xi _+^au',V_0}_{ \ol H^{\frac12-a+\varepsilon }(\rnp) ,\dot H^{a-\frac12-\varepsilon
}(\crnp)}\\
&=\ang{w,\Xi _+^aV_0}_{ \ol H^{\frac12+\varepsilon }(\rnp) ,\dot H^{-\frac12-\varepsilon
}(\crnp)}=\ang{w,\psi _0(x')\otimes \delta (x_n)}_{ \ol H^{\frac12+\varepsilon }(\rnp) ,\dot H^{-\frac12-\varepsilon
}(\crnp)}.
\endaligned\tag3.18
$$

Recall moreover from distribution
theory (cf.\ e.g.\ \cite{G09} p.\ 307) that the ``two-sided'' trace operator $\widetilde\gamma _0\colon
v(x)\mapsto \widetilde \gamma _0 v=v(x',0)$ has the mapping
$\widetilde\gamma _0^*\colon \varphi (x')\mapsto \varphi (x')\otimes
\delta (x_n)$ as adjoint, with continuity properties
$$
\widetilde \gamma _0\colon H^{\frac12+\varepsilon }({\Bbb R}^n)\to H^{\varepsilon }({\Bbb
R}^{n-1}),\quad
\widetilde \gamma _0^*\colon H^{-\varepsilon }({\Bbb R}^{n-1})\to H^{-\frac12-\varepsilon }({\Bbb
R}^{n}), \text{ for }\varepsilon >0. %\tag3.6
$$
Here $\widetilde\gamma _0^*\varphi $ is supported in $\{x_n=0\}$,
 hence lies in $\dot H^{-\frac12-\varepsilon }(\crnp)$.
Since $w\in \ol H^{\frac12+\varepsilon }(\rnp)$, it has an extension
$W\in H^{\frac12+\varepsilon }(\R^n)$ with
$w=r^+W$, and  $\gamma _0 w= \widetilde
\gamma _0W$.
Then
$$
\aligned
I&=\ang{w,\psi _0(x')\otimes \delta (x_n)}_{ \ol H^{\frac12+\varepsilon }(\rnp) ,\dot H^{-\frac12-\varepsilon
}(\crnp)}=\ang{W,\psi _0(x')\otimes \delta (x_n)}_{  H^{\frac12+\varepsilon }(\Rn) , H^{-\frac12-\varepsilon
}(\Rn)}\\
&=\ang{W,\widetilde \gamma _0^*\psi _0}_{  H^{\frac12+\varepsilon }(\Rn) , H^{-\frac12-\varepsilon
}(\Rn)}=\ang{\widetilde\gamma _0W,\psi _0}_{H^\varepsilon
(R^{n-1}),H^{-\varepsilon }(R^{n-1})}
=(\gamma
_0w,v _0)_{L_2(\R^{n-1})},
\endaligned
$$
since $\psi _0=v_0\in L_2({\Bbb R}^{n-1})$.
Finally, since $u'=U_1+u''$ with  $\gamma _0w=\gamma _0\Xi
_+^aU_1=\varphi _1=u_1+a\ang{D'}u_0$ by (3.6), (3.7), (3.13), we conclude:

\proclaim{Lemma 3.2} Let $P=(1-\Delta )^a$, and let $u,v\in
H^{(a-1)(s)}(\crnp)$, $s>a+\frac12$, 
%expanded as in {\rm (3.1)} 
%(3.6), (3.13)} 
%with $k=1$. 
with Dirichlet and Neumann boundary values $\gamma _0^{a-1}u=u_0$ and
$\gamma _1^{a-1}u=u_1$ (and similarly for $v$), as defined above. 
Let $V_0=\Cal F^{-1}_{\xi '\to x'}[\hat v_0(\xi ')I^{a-1}(x_n)e^{-\ang{\xi
'}x_n}]$. Then
$\int_{\rnp}r^+Pu\,\bar V_0\,dx$, understood as the duality {\rm (3.17)},
satisfies
$$
\aligned
\ang{r^+Pu,V_0}_{ \ol H^{\frac12-a+\varepsilon }(\rnp) ,\dot H^{a-\frac12-\varepsilon
}(\crnp)}&=\ang{r^+Pu',V_0}_{ \ol H^{\frac12-a+\varepsilon }(\rnp) ,\dot H^{a-\frac12-\varepsilon
}(\crnp)}\\
%&=(\varphi _1
%,v _0)_{L_2(\R^{n-1})}
&=(u_1+a\ang{D'}u_0
,v _0)_{L_2(\R^{n-1})}.
\endaligned\tag3.19
$$

\endproclaim

From this we obtain the Green's formula:

\proclaim{Theorem 3.3} Let $P=(1-\Delta )^a$, and let $u,v\in  H^{(a-1)(s)}(\crnp)$ with
$s>a+\frac12$.
Then 
$$
\multline
\ang{r^+Pu,v}_{ \ol H^{\frac12-a+\varepsilon }(\rnp) ,\dot H^{a-\frac12-\varepsilon
}(\crnp)}-\ang{u,r^+Pv}_{ \dot H^{a-\frac12-\varepsilon
}(\crnp),\ol H^{\frac12-a+\varepsilon }(\rnp) }\\
=\int_{{\Bbb R}^{n-1}}(\gamma _1^{a-1}u\,\gamma _0^{a-1}\bar
v-\gamma _0^{a-1}u\,\gamma _1^{a-1}\bar v)\, dx'.
\endmultline\tag3.20
$$
Here when $s\ge 2a$, the left-hand side can be written as an ordinary
integral
$$
\int_{\rnp}(r^+Pu\bar v\, dx -u\,\overline{r^+P v})\,
dx.
\tag3.21
$$

\endproclaim

\demo{Proof} It follows from (3.15) and Lemma 3.2 together  that when $s>a+\frac12$,
$$
\multline
\ang{r^+Pu,v}_{ \ol H^{\frac12-a+\varepsilon }(\rnp) ,\dot H^{a-\frac12-\varepsilon
}(\crnp)}\\
=\ang{r^+Pu',v'}_{ \ol H^{\frac12-a+\varepsilon }(\rnp) ,\dot H^{a-\frac12-\varepsilon
}(\crnp)}+(u_1 +a\ang{D'}u_0
,v _0)_{L_2(\R^{n-1})}
.\endmultline\tag3.22
$$
There is a similar formula obtained by interchanging $u$ and $v$ and conjugating:
$$
\multline
\ang{u,r^+Pv}_{ \dot H^{a-\frac12-\varepsilon
}(\crnp),\ol H^{\frac12-a+\varepsilon }(\rnp) }\\
=\ang{u',r^+Pv'}_{ \dot H^{a-\frac12-\varepsilon
}(\crnp),\ol H^{\frac12-a+\varepsilon }(\rnp) }
+(u_0
,v _1+a\ang{D'}v_0)_{L_2(\R^{n-1})}.
\endmultline
\tag3.23
$$
By Th.\ 4.1 of \cite{G16},
$$
\ang{r^+Pu',v'}_{ \ol H^{\frac12-a+\varepsilon }(\rnp) ,\dot H^{a-\frac12-\varepsilon
}(\crnp)}-\ang{u',r^+Pv'}_{ \dot H^{a-\frac12-\varepsilon
}(\crnp),\ol H^{\frac12-a+\varepsilon }(\rnp) }=0.
$$ 
 The dualities written there are consistent with the present ones, since  $u',v'\in H^{a(s)}(\crnp)=\Xi _+^{-a}e^+\ol
H^{\frac12+\varepsilon }(\rnp)\subset \dot H^{1-\varepsilon }$ and
$Pu', Pv'\in \ol H^{s-2a}(\rnp)\subset \ol H^{\frac12-a+\varepsilon }(\rnp)$.

Formula (3.20) then follows by taking the difference of (3.22) and
(3.23), using that $a\ang{D'}$ is selfadjoint.

If $s\ge 2a$, then $r^+Pu, r^+Pv\in \ol H^{s-2a}(\rnp)\subset
 L_2(\rnp)$.
Moreover, $u,v\in x_n^{a-1}\ol
H^{a+1}(\rnp)$ \linebreak $+\dot H^{2a}(\crnp)$,
cf.\ (2.12).
So $r^+Pu\,\bar v$ and $u\,\overline{r^+Pv}$ are functions, and
we can write the dualities as in (3.21), keeping the interpretation in
 mind.\qed
\enddemo

\example{Remark 3.4} We take the opportunity to mention that the formula
(5.14) in \cite{G15} Th.\ 5.4 is only exact when $M=1$; when $M>1$,
there 
are missing some terms with $x_n^{j+\mu }K_0\gamma _{\mu ,k}$ ($k<j$) and
$\psi $do coefficients, like $a\ang{D'}$ in (3.7) here. The conclusion (5.15) remains valid.
A corrected formula will be included in a forthcoming paper \cite{G18}.
\endexample

\subhead 4. Green's formula for variable-coefficient operators \endsubhead

Let $P=\OP (p(x,\xi ))$ be a classical $\psi $do on ${\Bbb R}^n$ of
order $2a>0$ (global estimates), with
symbol $p(x,\xi )\sim\sum_{j\in{\Bbb N}_0}p_j(x,\xi )$.
%, i.e., the
%$p_j$ are homogeneous in $\xi $ of degree $2a-j$ for $|\xi |\ge 1$, and
%there are global estimates 
We assume that $p$ is
{\it even}, i.e.,
$$
p_j(x,-\xi )=(-1)^jp_j(x,\xi )\text{ for }|\xi |\ge 1, \text{ all }x.%\tag5.1
$$
Let $\Omega $ be a smooth bounded subset of ${\Bbb R}^n$, or $\Omega =\rnp$. 
The evenness implies that $p$ satisfies (2.4) with $m=2a$, $\mu =a$,
so $p$ (or $P$) has the $a$-transmission
property at $\Omega $. The adjoint $P^*$ is likewise even.
(Evenness is assumed for simplicity in the formulations; everything
 goes through when $\Omega$ is given on beforehand and $P$ is just
 assumed to have the $a$-transmission property with respect to the
 particular $\Omega $.)

Green's formula for these operators %$x$-dependent
will first be shown in the case $\Omega =\rnp$, and
afterwards generalized to the curved case. The main strategy is to
reduce as much as possible to rules from the Boutet de Monvel calculus, where the
issues of operators passing to and from the boundary are dealt with in a
systematic way (in the 0-transmission case). We refer
the reader to e.g.\  \cite{G09} for a general presentation of the
calculus; a few important ingredients are collected in the Appendix
here.

\proclaim{Theorem 4.1} Let $P$ be a classical
$\psi $do on ${\Bbb R}^n$ of order $2a>0$ with even symbol. The following Green's
formula holds for  $u,v\in H^{(a-1)(s)}(\crnp)$ when $s>a+\frac12$:
$$\aligned
\ang{r^+Pu,v}&_{\ol H^{-a+\frac12+\varepsilon }(\rnp), \dot
H^{a-\frac12-\varepsilon }(\crnp)}-\ang{u,r^+P^*v}_{\dot H^{a-\frac12-\varepsilon }(\crnp),\ol
H^{-a+\frac12+\varepsilon }(\rnp)}\\
&=( s_0\gamma _1^{a-1}u,\gamma
_0^{a-1}v)-( s_0\gamma _0^{a-1}u,\gamma _1^{a-1}v)+(B\gamma
_0^{a-1}u,\gamma _0^{a-1}v),
\endaligned\tag4.1
$$
with $L_2({\Bbb R}^{n-1})$-scalar products in the right-hand side;
here $s_0=p_0(x',0,0,1)$, and $B$ is a first-order $\psi $do  on ${\Bbb
R}^{n-1}$
whose symbol equals the jump at $z_n=0$ in the bounded part of $\Cal
F^{-1}_{\xi _n\to z_n}q(x',0,\xi )$, where $q$ is the symbol of $\Xi _-^{-a}P\Xi _+^{-a}$.
The right-hand side can also be written
$$
\aligned
\Gamma (a)\Gamma (a+1)\int_{{\Bbb R}^{n-1}}(s_0\gamma
_1(\tfrac u{{x_n}^{a-1}})\gamma _0(\tfrac{\bar v}{{x_n}^{a-1}})
& -s_0\gamma _0(\tfrac{u}{{x_n}^{a-1}})\gamma
_1(\tfrac{\bar v}{{x_n}^{a-1}})+\\
&+aB\gamma _0(\tfrac{u}{{x_n}^{a-1}})\gamma _0(\tfrac{\bar v}{{x_n}^{a-1}}))
\, dx' .
\endaligned
\tag4.2
$$

When $s\ge 2a$, the dualities in the left-hand side can be written as
integrals over $\rnp $.
\endproclaim

\demo{Proof}
Define $Q=\Xi _-^{-a}P\Xi _+^{-a}$, it is a generalized $\psi $do of
order 0, with a symbol $q(x,\xi )$ that is the sum of
$s_0(x)=p_0(x,0,1)$ and a function in $S^0(\Cal H_{-1})$ (notation explained in
the Appendix). The principal symbol 
$q_0(x,\xi
)$ is is the $\psi $do symbol  $p_0(x,\xi )[\xi ]^{-2a}$. 
We also need the $(x',y_n)$-form $q'(x',y_n,\xi )$ related to $q$ by (A.11).
Now we can write $P$ as $P=\Xi _-^aQ \Xi
_+^a$.

(In the following study, $\Xi  _\pm^t$ and $\Lambda  _\pm^t$ can be used
equally well. With the use of $\Lambda _\pm^t$ and the corresponding
choice of $Q$, the calculations stay within
true pseudodifferential operators as much as possible. With
$\Xi_\pm^t$ the formulas are simpler and more direct; here we draw on
the fact that  
 $Q$ has symbol in $S^0(\Cal H_{-1})$ plus
smooth functions, where the rules for Poisson and trace operators we need
are still valid.) 

Using the description of $\gamma _0^{a-1}$ and $K_0^{a-1}$ given in
Section 3, we decompose
a function $u\in H^{(a-1)(s)}(\crnp)$ as 
$$
u=u'+K_0^{a-1}\gamma _0^{a-1}u;\text{ here } u'\in
H^{a(s)}(\crnp)\text{ since }\gamma _0^{a-1}u'=0.\tag4.3
$$
There is a similar decomposition for $v$, and we denote 
$$
\gamma _0^{a-1}u= u_0 ,\quad \gamma _0^{a-1}v= v_0 ,\text{ both in }
 H^{s-a+\frac12}({\Bbb R}^{n-1}).\tag4.4
$$
The idea of the proof is to show that the contributions from the terms
$K^a_0u_0$ and $K^a_0v_0$ give expressions that can be reduced to give
the right-hand side of (4.1). Here we eliminate the
fractional-order factors so that rules from the Boutet de Monvel
calculus can be applied.

We assume $s>a+\frac12$, so instead of $s$ we can insert
$a+\frac12+\varepsilon $, and $$\aligned
u,v,K_0^{a-1}u _0 \text{ and }K_0^{a-1}v _0 
&\in 
H^{(a-1)(a+\frac12+\varepsilon )}(\crnp)\subset \dot H^{a-\frac12-\varepsilon }(\crnp),\\
 u_0,v_0 &\in H^{1+\varepsilon }({\Bbb R}^{n-1}),\\
u',v'&\in H^{a(a+\frac12+\varepsilon )}(\crnp)\subset \dot H^{a+\frac12-\varepsilon }(\crnp),\\
r^+Pu, r^+P^*v&\in \ol H^{-a+\frac12+\varepsilon }(\rnp).
\endaligned\tag4.5$$
Then $\ang{r^+Pu,v}$ can be interpreted as
$$
\ang{r^+Pu,v}=\ang{r^+Pu,v}_{\ol H^{-a+\frac12+\varepsilon }, \dot
H^{a-\frac12-\varepsilon }}.\tag4.6
$$
This expression is split into four parts by applying (4.3) to $u$ and $v$:
$$
\aligned
\ang{r^+Pu,v}&=I_1+I_2+I_3+I_4, \\
I_1&=\ang{r^+Pu',v'},\quad I_2=\ang{r^+PK_0^{a-1}u_0
,v'},\\
I_3&=\ang{r^+Pu',K_0^{a-1} v_0  },\quad I_4=\ang{r^+PK_0^{a-1} u_0 ,K_0^{a-1} v_0  }.
\endaligned\tag4.7
$$

$I_1$ will be kept unchanged, to match a similar term with $P^*$ later.

For $I_2$ we observe (for small $\varepsilon $), using (2.8):
$$
\aligned
I_2&=\ang{r^+PK_0^{a-1} u_0 ,v'}_{\ol H^{-a+\frac12+\varepsilon }, \dot
H^{a-\frac12-\varepsilon }}\\
&=\ang{r^+\Xi_-^{a}e^+r^+Q\Xi_+^{a}\Xi_+^{1-a}e^+K_{0}u_0
,v'}_{\ol H^{-a+\frac12+\varepsilon }, \dot H^{a-\frac12-\varepsilon
}}\\
&=\ang{r^+Q\Xi_+^{a}\Xi_+^{1-a}e^+K_{0} u_0 ,\Xi 
_+^av'}_{\ol H^{\frac12+\varepsilon }, \dot H^{-\frac12-\varepsilon
}}\\
&=\ang{r^+Q\Xi_+^{1}e^+K_{0} u_0 ,\Xi 
_+^av'}_{\ol H^{\frac12+\varepsilon }, \dot H^{-\frac12-\varepsilon }}.
\endaligned\tag4.8
$$

We shall now apply the rules of calculus for $\psi $dbo's, as recalled
in the Appendix. Let us mention here that the projection (idempotent) $h^+$ applied
to $\xi _n$-dependent symbols can just be thought of as the Fourier
transform of the projection  $e^+r^+$ in $L_2({\Bbb R})$. It is
applied to a more refined space $\Cal H=\Cal F_{x_n\to \xi _n}(e^-\Cal
S_-\oplus e^+\Cal S_+\oplus {\Bbb C}[\delta ])$, where $\Cal
S_\pm$ is short for $r^\pm\Cal S({\Bbb R}) $, and ${\Bbb C}[\delta ]$
is the space of distributions supported in $\{x_n=0\}$. Then $\Cal
H=\Cal H^+\oplus\Cal H^-$ where $\Cal H^+=\Cal F _{x_n\to \xi
_n}(e^+\Cal S_+)$ and $\Cal H^-=\Cal F _{x_n\to \xi _n}(e^-\Cal
S_-)\oplus {\Bbb C}[\xi _n]$, with ${\Bbb C}[\xi _n]$ denoting the space of
complex polynomials in $\xi _n$. Now $h^\pm$ is the projection of
$\Cal H$ onto $\Cal H^\pm$ along $\Cal H^\mp$. More details in the Appendix, and a
full deduction e.g.\ in \cite{G09} Sect.\ 10.2.

Using the symbol $q'$ of $Q$ in $(x',y_n)$-form, we have: 
$$
r^+Q\Xi_+^{1}e^+K_{0} =\operatorname{OPK}(h^+(q'\#\chi _+^1\chi
_+^{-1}))=\operatorname{OPK}(h^+q'(x',0,\xi )),\tag4.9
$$
which leads to
$$
I_2=
\ang{\operatorname{OPK}(h^+q')u_0, \Xi 
_+^av'}_{\ol H^{\frac12+\varepsilon }, \dot H^{-\frac12-\varepsilon }}.\tag4.10
$$
If $Q=I$ we get zero here, but for general $Q$ there is a
 nontrivial contribution from $h^+q'$. Note that with $s_0(x)=p_0(x,0,1)$,
$$
q(x,\xi )=s_0(x)+h_{-1}q(x,\xi ),\quad h^+q=h^+h_{-1}q,\quad h^-q=s_0+h^-_{-1}q,\tag4.11
$$
with similar rules for $q'(x',y_n,\xi )$.

Now consider $I_3$. Here, using (2.8) and the adjoint $K_0^*$ of
$K_0$,
$$
\aligned
I_3&=\ang{r^+Pu',K_0^{a-1} v_0  }=\ang{r^+\Xi_-^aQ\Xi_+^au',\Xi_+^{1-a}e^+K_{0} v_0  }_{\ol H^{-a+\frac12+\varepsilon }, \dot
H^{a-\frac12-\varepsilon }}\\
&=\ang{r^+\Xi_-^{1-a}e^+r^+\Xi _-^{a}e^+r^+Q\Xi_+^au',K_{0} v_0  }_{\ol H^{-\frac12+\varepsilon }, \dot
H^{\frac12-\varepsilon }}\\
&=\ang{r^+\Xi_-^1e^+r^+Q\Xi_+^au',K_{0} v_0  }_{\ol H^{-\frac12+\varepsilon }, \dot
H^{\frac12-\varepsilon }}\\
&=\ang{K_0^*r^+\Xi_-^1Q\Xi_+^au', v_0  }_{ H^{\varepsilon }({\Bbb R}^{n-1}), 
H^{-\varepsilon }({\Bbb R}^{n-1})}
=(K_0^*r^+\Xi_-^1Q\Xi_+^au', v_0  )_{L_2({\Bbb R}^{n-1})},
\endaligned\tag4.12$$
since $ v_0  $ is in $L_2({\Bbb R}^{n-1})$.
It is
used that $\ol H^t$ identifies with $\dot H^t$ for $|t|<\frac12$
(there the
indication $e^+$ is understood). 

Denote $\Xi_+^au'=w\in e^+\ol H^{\frac12+\varepsilon }(\rnp)$. 
Observe that
$$
K_0^*r^+\Xi_-^1Qw=\operatorname{OPT}(h^-(\chi _-^{-1}\chi
_-^{1}\#q))w=
\operatorname{OPT}(h^-q(x',0,\xi ))w,
$$
by the rules of calculus, so 
$$
\aligned
I_3&=(\operatorname{OPT}(h^-q)\Xi_+^au', v_0 
)_{L_2({\Bbb R}^{n-1})}=((s_0\gamma _0+\operatorname{OPT}(h_{-1}^-q))\Xi_+^au', v_0 )_{L_2({\Bbb R}^{n-1})}\\
&=
(s_0\gamma _0^au', v_0 )_{L_2({\Bbb R}^{n-1})}+(\operatorname{OPT}(h^-_{-1}q)\Xi_+^au', v_0  )_{L_2({\Bbb R}^{n-1})}.
\endaligned\tag4.13
$$
The first term  is expected from Theorem 3.3, and there is a nontrivial
extra term.

Finally, consider $I_4$: Here
$$
\aligned
 I_4&=\ang{r^+PK_0^{a-1} u_0 ,K_0^{a-1} v_0  }_{\ol H^{-a+\frac12+\varepsilon }, \dot
H^{a-\frac12-\varepsilon }}\\
&=\ang{r^+\Xi 
 _-^aQ\Xi_+^a\Xi_+^{1-a}e^+K_0 u_0 ,\Xi 
 _+^{1-a}e^+K_0 v_0  }_{\ol H^{-a+\frac12+\varepsilon }, \dot
H^{a-\frac12-\varepsilon }}\\
&=\ang{r^+\Xi_-^1e^+r^+Q\Xi_+^1e^+K_0 u_0 ,K_0 v_0  }_{\ol H^{-\frac12+\varepsilon }, \dot
H^{\frac12-\varepsilon }}\\
%&=\ang{r^+Q\Xi_+^1e^+K_0 u_0 ,\Xi _+^1K_0 v_0  }_{\ol H^{\frac12+\varepsilon },% \dot
%H^{-\frac12-\varepsilon }} [nyt]\\
%&\text{passer det med BdM teorien? vi kan i stedet regne p\aa{} }S_0\\
&=\ang{K_0^*r^+\Xi_-^1Q\Xi_+^1e^+K_0 u_0 , v_0  }_{ H^{\varepsilon }({\Bbb R}^{n-1}), 
H^{-\varepsilon }({\Bbb R}^{n-1})}\\
&=(\Cal B u_0 , v_0  )_{L_2({\Bbb R}^{n-1})},
\endaligned\tag4.14
$$
where ${\Cal B}=K_0^*r^+\Xi_{-}^1e^+r^+Q\Xi_+^1e^+K_0$ is a certain $\psi $do
on ${\Bbb R}^{n-1}$ of order 1. We can reduce this expression by rules
of calculus involving the so-called  plus-integral, cf.\ (A.14)ff. and (A.15). 
As in (4.9), 
the symbol of the Poisson operator
$r^+Q\Xi _+^{1}e^+K_0$ is $h^+q'(x',0,\xi )$, which by composition 
with $r^+\Xi_{-}^1e^+$ to the left gives a Poisson operator with
symbol $h^+(\chi _-^1\#h^+q')$; hence
the symbol
$\frak b(x',\xi ')$ of ${\Cal B}$ satisfies by (A.14),
$$
\aligned
\frak b(x',\xi ')&=\tfrac1{2\pi }\int^+\chi _-^{-1}\#h^+(\chi _-^1\#h^+q'(x',0,\xi ))\,
d\xi _n\\
&=\tfrac1{2\pi }\int^+\chi _-^{-1}\#(\chi _-^1\#h^+q'(x',0,\xi )-h^- (\chi _-^1\#h^+(q'(x',0,\xi )))\,
d\xi _n\\
&=\tfrac1{2\pi }\int^+(h^+q'(x',0,\xi )-\chi _-^{-1}\#h^- (\chi _-^1\#h^+(q'(x',0,\xi )))\,
d\xi _n\\
&=\tfrac1{2\pi }\int^+ h^+q'(x',0,\xi )\,
d\xi _n;
\endaligned\tag4.15
$$
it was used here %moreover 
that the plus-integral vanishes on $\Cal H^-$.
(The Leibniz product $\#$ pertains to the $x'$-dependence.)
Now observe that with 
 $\Cal F^{-1}_{\xi _n\to z_n}q'(x',0,\xi )$ denoted $\check {q'}(x',0,\xi ',z_n)$, 
$$
\tfrac1{2\pi }\int^+ h^+q'(x',0,\xi )\,
d\xi _n =\lim_{z_n\to
0+}\check {q'}(x',0,\xi ',z_n),
$$
cf.\ (A.15). (We use here for each fixed $(x',\xi ')$ that since $q'(x',0,\xi ',\xi _n)$
is the sum of  a function $f(x',\xi ',\xi _n)\in
\Cal H_{-1}$ and the constant $s_0(x',0)$,
$\Cal
F^{-1}_{\xi _n\to z_n}q'(x',0,\xi ',\xi _n)$ is the sum of a function 
$\check f(x',\xi ',z_n)\in e^-\Cal S(\crm)\oplus e^+\Cal S(\crp)$ and the distribution $s_0(x',0)\delta (z_n)$, where
the latter does not enter in the limit from the right.) 
Moreover, in view of the formula (A.11),
$$
\check {q'}(x',0,\xi ',z_n)\sim \sum_{j\in {\Bbb N}_0}\tfrac
1{j!}z_n^j\partial_{x_n}^j\check  q(x',0,\xi ',z_n),
$$
hence
$$
\lim_{z_n\to
0+}\check {q'}(x',0,\xi ',z_n)=\lim_{z_n\to
0+}\check {q}(x',0,\xi ',z_n),
$$
since the positive powers of $z_n$ vanish at 0. It follows that
$$
\frak b(x',\xi ')=\lim_{z_n\to
0+}\check {q}(x',0,\xi ',z_n)=\tfrac1{2\pi }\int^+ h^+q(x',0,\xi )\,
d\xi _n .\tag4.16
$$

There is a similar decomposition of $\ang{u,r^+P^*v}$ in four terms 
$I'_1,I'_2,I'_3,I'_4$. We here note that $Q^*$ has the symbol in
$y$-form $\bar q(y,\xi )$, where
$$
\bar q=\bar s_0+\overline{h_{-1}q},\quad  h^+\bar
q=h^+(\overline{h_{-1}q})=\overline {h^-_{-1}q},
\quad  h^-\bar q=\bar s_0+h^-(\overline{h_{-1}q})=\bar s_0+\overline{h^+q}.\tag4.17
$$
The symbol of $Q^*$ in $(y',x_n)$-form is $\bar q'(y',x_n,\xi )$,
satisfying similar rules.

We here find the formulas
$$
\aligned
I'_1&=\ang{u',r^+P^*v'}_{\dot
H^{a-\frac12-\varepsilon },\ol H^{-a+\frac12+\varepsilon }},\\
I'_2&=\ang{\Xi 
_+^au',\operatorname{OPK}(\overline{h^-_{-1}q(y',0,\xi )}) v_0 , }_{\dot
H^{-\frac12-\varepsilon },\ol H^{\frac12+\varepsilon }},\\
I_3'&=\ang{ u_0 ,\bar s_0\gamma _0^av'}
+\ang{ v_0  ,\operatorname{OPT}(\overline{h^+q'(y',0,\xi )})\Xi_+^av' },\\
I_4'&=( u_0 ,{\Cal B}' v_0  )_{L_2({\Bbb R}^{n-1})},
\endaligned\tag4.18$$
where the operators derived from $q$ are in $y'$-form.
The operator ${\Cal B}'$ is defined as
$$
{\Cal B}'=K_0^*r^+\Xi_{-}^1e^+r^+Q^*\Xi_+^1e^+K_0;
$$
it has the symbol (reduced as in (4.15))
$$
\aligned
\frak b'(y',\xi ')&=\tfrac1{2\pi }\int^+\chi _-^{-1}h^+(\chi _-^1h^+\bar q(y',0,\xi ))\,
d\xi _n\\
&=\tfrac1{2\pi }\int^+h^+\bar q(y',0,\xi )\,
d\xi _n = \lim_{z_n\to
0+}\Cal F^{-1}_{\xi _n\to z_n}\bar q(y',0,\xi ',z_n).
\endaligned%\tag4.19
$$
In view of (4.17), $h^+\bar q= \overline {h^-_{-1}q} $. We note that
$h^+q$ and $\overline{h^-_{-1}q}$ can be quite different, so ${\Cal B}$
and ${\Cal B}'$ (or its adjoint) are in general different from one another. 

The adjoint of ${\Cal B}'$ is ${\Cal B}^{\prime
*}=\OP'(\overline{\frak b'}(x',\xi
'))$. 
%Denote the reflection operator by $J$, $J\colon \varphi
%(t)\mapsto \varphi (-t)$. 
It is well-known (and easily checked, cf.\ e.g.\ \cite{G09} p.\ 118) that
for a distribution $\varphi (\xi _n)$,  $\overline{[\Cal F^{-1}_{\xi
_n\to z_n}\overline\varphi](z_n)} =[\Cal F^{-1}_{\xi _n\to z_n}\varphi ](-z_n)$. Then
$$
\aligned
I'_4&=({\Cal B}^{\prime*}u_0,v_0),\text{ where }{\Cal B}^{\prime
*}=\OP'(\overline{\frak b'}(x',\xi
')),\\
\overline{\frak b'}(x',\xi ')&=
\lim_{z_n\to
0+}\overline{\Cal F^{-1}_{\xi _n\to z_n}\bar q(x',0,\xi ',z_n)}=
\lim_{z_n\to
0+}\check q(x',0,\xi ',-z_n)\\
&=\lim_{z_n\to
0-}\check q(x',0,\xi ',z_n).
\endaligned\tag4.19
$$

The full right-hand side in (4.1) is
$I_1+I_2+I_3+I_4-I'_1-I'_2-I'_3-I'_4$, that we can now calculate. 
Here  we find:
$$
I_1-I'_1=0,
$$
by \cite{G16} Th.\ 4.1. Next,
$$
\aligned
I_2-I'_2&=\ang{\operatorname{OPK}(h^+q')u_0, \Xi
_+^av'}_{\ol H^{\frac12+\varepsilon }, \dot H^{-\frac12-\varepsilon
}}-\ang{\Xi_+^au',\operatorname{OPK}(\overline{h^-_{-1}q}) v_0 }_{\dot
H^{-\frac12-\varepsilon },\ol H^{\frac12+\varepsilon }}\\
&=\ang{\operatorname{OPK}(h^+q')u_0, \Xi
_+^av'}_{\ol H^{\frac12+\varepsilon }, \dot H^{-\frac12-\varepsilon
}}-\ang{\operatorname{OPT}({h^-_{-1}q})\Xi_+^au', v_0 }_{
H^{-\varepsilon }({\Bbb R}^{n-1}), H^{\varepsilon }({\Bbb R}^{n-1})}
, \\
I_3-I_3'&=\ang{s_0\gamma _0^au', v_0 
}+\ang{\operatorname{OPT}(h^-_{-1}q)\Xi_+^au', v_0  }-\ang{ v_0  ,\bar s_0\gamma
_0^av'  }-\ang{ u_0 ,\operatorname{OPT}(\overline{h^+q'})\Xi_+^av' }\\
&=\ang{s_0(\gamma _1^{a-1}u+a\ang{D '}u_0),v_0
}-\ang{s_0u_0 ,\gamma
_1^{a-1}v +a\ang{D'}v_0 }\\
&\quad+\ang{\operatorname{OPT}(h^-_{-1}q)\Xi_+^au', v_0  }
-\ang{ u_0 ,\operatorname{OPT}(\overline{h^+q'})\Xi_+^av' }\\
&=\ang{s_0\gamma _1^{a-1}u,v_0
}-\ang{s_0u_0 ,\gamma
_1^{a-1}v  }\\
&\quad+\ang{\operatorname{OPT}(h^-_{-1}q)\Xi_+^au', v_0  }
-\ang{\operatorname{OPK}({h^+q'}) u_0 ,\Xi_+^av' };
\endaligned
$$
we have here used the rules for adjoints of Poisson and trace
operators,
and the fact that $\gamma _0^au'=\gamma _1^{a-1}u+a\ang{D'}u_0$, cf.\ (3.9). (Operators in $y'$-form in the right-hand side give operators in $x'$-form when transposed to the
left-hand side.) Thus
$$
%\aligned
I_2+I_3-I'_2-I'_3=\ang{s_0\gamma _1^{a-1}u,\gamma _0^{a-1}v
}-\ang{s_0\gamma _0^{a-1}u ,\gamma
_1^{a-1}v  }.
%\endaligned
$$

Finally,
$$
I_4-I'_4=(({\Cal B}-{{\Cal B}'}^*)u_0,v_0)_{L_2({\Bbb R}^{n-1})},=(({\Cal B}-{{\Cal B}'}^*)\gamma _0^{a-1}u,\gamma _0^{a-1}v)_{L_2({\Bbb R}^{n-1})},
$$
where  $B={\Cal B}-{{\Cal B}'}^*$ satisfies, with $b(x',\xi ')=\frak
b(x',\xi ')-\overline{\frak b'}(x',\xi ')$,
$$
B=\OP'(b(x',\xi ')),\quad b(x',\xi ')=\lim_{z_n\to
0+}\check q(x',0,\xi ',z_n)-\lim_{z_n\to
0-}\check q(x',0,\xi ',z_n),\tag4.20
$$
the jump at $z_n=0$ in the bounded part of $\check q(x',0,\xi ',z_n)$.
Altogether, we find:
$$
\ang{r^+Pu,v}-\ang{u,r^+{P^*v}}
=( s_0\gamma _1^{a-1}u,\gamma
_0^{a-1}v)-( s_0\gamma _0^{a-1}u,\gamma _1^{a-1}v)+(B\gamma
_0^{a-1}u,\gamma _0^{a-1}v),
$$
with $L_2({\Bbb R}^{n-1})$-dualities in the right-hand side,
containing a generally nontrivial $\psi $do $B$ on ${\Bbb R}^{n-1}$. 
The formulation in (4.2) follows by inserting (3.3).

The last assertion is seen as in Theorem 3.3. \qed
\enddemo

\example{Remark 4.2} Theorem 3.3 shows that $B=0$ in the simple case
 of $P=(1-\Delta )^a$. More generally, for operators with principal
 symbol $|\xi |^{2a}$, the principal (first-order) symbol of $B$ is zero,
so $B$ is of order 0. 
 
Note moreover that if $q(x',0,\xi )-s_0(x',0)$ is $O(\ang{\xi _n}^{-2})$, hence
integrable in $\xi _n$, for all $\xi '$, then $\check q(x',0,\xi
',z_n)-s_0(x',0)\delta (z_n)$
is a continuous function of $z_n$, so $b(x',\xi ')=0$ and hence $B=0$.
In general, $q$ is the sum of a part depending only on the principal and
subprincipal terms in $p$ and a part that is  $O(\ang \xi ^{-2})$; 
the latter part does not contribute to $B$.
\endexample

To extend the result of Theorem 4.1 to a domain $\Omega $ with curved boundary, we
shall use a suitable cover of $\comega$ by coordinate charts, and a
suitable partition of unity. Such choices were described in
\cite{G16}, and we recall them in the following remark.

\example{Remark 4.3}
$\comega$ has a finite cover by bounded open sets $U_0,\dots, U_{I}$ with
$C^\infty $-dif\-fe\-o\-mor\-phisms
$\kappa _i\colon U_i\to V_i$, $V_i$ bounded open in ${\Bbb R}^n$, such
 that $U_i^+=U_i\cap \Omega $ is mapped to $V_i^+=V_i\cap\rnp$ and
 $U_i'=U_i\cap\partial\Omega $ is mapped to $V_i'=V_i\cap\partial\crnp$; as
 usual we write $\partial\crnp={\Bbb R}^{n-1}$. 
For any such cover there exists an associated partition of unity, namely a family
of functions $\varrho _i\in C_0^\infty (U_i)$ taking values in $[0,1]$
such that $\sum_{i=0,\dots,I}\varrho _i$ is 1 on a neighborhood of $\comega$.

When $P$ is a $\psi
 $do on ${\Bbb R}^n$, its application to functions
 supported in $U_i$ carries over to functions on $V_i$ as a $\psi $do
 $\underline P^{(i)}$ defined by 
$$
\underline P^{(i)}v=P(v\circ \kappa _i)\circ \kappa _i^{-1},\quad v\in C_0^\infty (V_i).\tag4.21
$$

We shall use a convenient system of coordinate charts as described 
in \cite{G16}, Remark
4.3: Here $\partial\Omega $ is covered with coordinate charts  $\kappa '_i\colon U'_i\to V'_i\subset {\Bbb R}^{n-1}$,
$i=1,\dots, I_0$, and the $\kappa _i$ will be defined on certain subsets  of a tubular
neighborhood $\Sigma _r=\{x'+t\nu (x')\mid x'\in\partial\Omega , |t|<r\}$,
where $\nu (x')=(\nu _1(x'),\dots,\nu _n(x'))$ is the interior normal
to $\partial\Omega $ at $x'\in\partial\Omega $, and $r$ is taken so
small that the mapping $ x'+t\nu (x')\mapsto (x',t)$ is a
diffeomorphism from $\Sigma _r$ to $\partial\Omega \times
\,]-r,r[\,$. For each $i$, $\kappa _i$ is defined as the mapping $\kappa
_i\colon x'+t\nu (x')\mapsto (\kappa '_i(x'),t)$ ($x'\in U'_i$).  $\kappa _i$ goes from $U_i$ to
$V_i$, where
$$
U_i=\{x'+t\nu (x')\mid x'\in U'_i, |t|<r\}, \quad V_i=V'_i\times
\,]-r,r[\,.\tag4.22
$$
These charts are supplied with a chart
consisting of the identity mapping on an open set $U_0$ containing $\Omega
\setminus \overline \Sigma _{r,+}$, with $\overline U_0\subset \Omega
$, to get a full cover of $\comega$.

Note that the normal $\nu (x')$ at $x
'\in\partial\Omega $ is carried over to the normal $(0,1)$ at $(\kappa
'_i(x'),0)$ when $x'\in U_i'$. The halfline $L_{x'}=\{x'+t\nu (x')\mid t\ge 0\}$ is the
geodesic into $\Omega $ orthogonal to $\partial\Omega $ at $x'$ (with
respect to the Euclidean metric on ${\Bbb R}^n$), and there is a positive
$r'\le r$ such that for $0<t<r'$, the distance $d(x)$ between  $x=x'+t\nu
(x')$ and $\partial\Omega $
equals $t$. Then $t$ plays the role of $d$ in the definition of
expansions and boundary values of $u\in \Cal E_{a-1}(\comega)$ in
\cite{G15} (5.3)ff.\ (cf.\ also (3.1) above) 
$$
u=\tfrac 1{\Gamma (a)}t^{a-1}u_0+\tfrac 1{\Gamma (a+1)}t^{a}u_1+\tfrac
1{\Gamma (a+2)}t^{a+1}u_2+\dots\text{ for }t>0,\quad u=0\text{ for }t<0,\tag4.23
$$
where the $u_j$ are constant in $t$ for $t<r'$;  this serves to define
the boundary values $\gamma ^{a-1}_{j}u=\gamma _0u_j\,(=u_j|_{t=0})$, $j=0,1,2,\dots$.
The definition extends to define the two first boundary values
$\gamma _0^{a-1}u$ and $\gamma _1^{a-1}u$ when $u\in
H^{(a-1)(s)}(\comega)$ with $s>a+\frac12$ (for the first boundary
value, $s>a-\frac12$ suffices). By comparison of (4.23) with
$t^{a-1}$ times the
Taylor expansion of $u/t^{a-1}$ in $t$, we also have:
$$
\gamma _0^{a-1}u=\Gamma (a)\gamma _0(u/t^{a-1}),\quad \gamma _1^{a-1}u=\Gamma (a+1)\gamma _1(u/t^{a-1})=\Gamma (a+1)\gamma _0(\partial_t(u/t^{a-1})),\tag4.24
$$
similarly as in (3.3).

In addition to the above construction of a cover by coordinate charts
and an associated partition of unity, it is for some purposes
practical to have a
 partition of unity {\it subordinate to a cover} as in
\cite{G16} Lemma 4.4. The cover is constructed
from the cover we have just described, by an
augmentation by extra
coordinate charts $\kappa _i\colon U_i\to V_i$, $i=I_0+1,\dots, I_1$,
such that there is a
par\-ti\-tion of unity $\varrho _k$, $k=1,\dots, J_0$,
where for any two functions $\varrho _k,\varrho _l$ there is an
$i=i(k,l)$ in $\{0,1,\dots,I_1\}$ for which $\supp \varrho _k\cup\supp \varrho
_l\subset U_{i(k,l)}$, see details in \cite{G16}. The maps $\kappa _i\colon U_i\to V_i$ still have the property that
the normal coordinate $t$ at the boundary goes over into $x_n$
(without distortion) for small $t$. 
\endexample

\proclaim{Theorem 4.4}  Let $P$ is a classical
$\psi $do on ${\Bbb R}^n$ of order $2a>0$ with even symbol, and let
$\Omega $ be a smooth bounded subset. The
following Green's formula holds for $u,v\in H^{(a-1)(s)}(\comega)$
when  $s>a+\frac12$:
$$\aligned
\ang{r^+Pu,v}&_{\ol H^{-a+\frac12+\varepsilon }(\Omega ), \dot
H^{a-\frac12-\varepsilon }(\comega)}-\ang{u,r^+P^*v}_{\dot H^{a-\frac12-\varepsilon }(\comega),\ol
H^{-a+\frac12+\varepsilon }(\Omega )}\\
&=( s_0\gamma _1^{a-1}u,\gamma
_0^{a-1}v)-( s_0\gamma _0^{a-1}u,\gamma _1^{a-1}v)+(B\gamma
_0^{a-1}u,\gamma _0^{a-1}v);
\endaligned\tag4.25
$$
the scalar products in the right-hand side are in 
$L_2(\partial\Omega )$, $s_0(x)=p_0(x,\nu (x))$ at boundary points $x$,
and $B$ is a first-order $\psi $do on $\partial\Omega $, depending
only on the principal and subprincipal symbols of $P$.

Here the right-hand side equals, in terms of the normal coordinate $t$ 
(with $\gamma _1w=\gamma _0\partial_tw$): 
$$
%\aligned
\Gamma (a)\Gamma (a+1)\int_{\partial\Omega }(s_0\gamma
_1(\tfrac u{t^{a-1}})\gamma _0(\tfrac{\bar v}{t^{a-1}})
 -s_0\gamma _0(\tfrac{u}{t^{a-1}})\gamma
_1(\tfrac{\bar v}{t^{a-1}})+aB\gamma _0(\tfrac{u}{t^{a-1}})\gamma _0(\tfrac{\bar v}{t^{a-1}}))
\, d\sigma .
%\endaligned
\tag4.26
$$

When $s\ge 2a$, the left-hand
side can be written as an integral over $\Omega $.

\endproclaim

\demo{Proof} For this proof we will use the cover $U_i$,
$i=0,\dots,I_1$, and the subordinate partition of
unity $\varrho _k$, $k=1,\dots,J_0$, described at the end of Remark
4.3. 

The first step is to show that the problem can be localized, i.e.\
that it suffices to prove the formula for functions supported in one
of the coordinate patches $U_i$. This is not completely obvious, since
the space $H^{(a-1)(s)}(\comega)$ is of the form $\Lambda _+^{-(a-1)}
e^+\ol H^{s-a+1}(\Omega )$, where $e^+\ol H^{s-a+1}(\Omega )$ is
preserved under multiplication by cutoff functions, but $\Lambda
_+^{-(a-1)}$ is nonlocal. We proceed as follows:

Choose nonnegative functions $\psi _k,\zeta  _k\in C_0^\infty
(U_i)$ such that
$\zeta  _k\varrho _k=\varrho _k$, i.e., $\zeta _k$ is 1 on $\supp \varrho
_k$, and similarly $\psi _k\zeta _k=\zeta _k$.
Let $u\in H_p^{(a-1)(s)}(\comega)$. Then $u=\Lambda _+^{(-a+1)}z$ for some
$z\in e^+\ol H_p^{s-a+1}(\Omega )$, and we can write
$$
\aligned
u&=\Lambda _+^{(-a+1)}{\sum}_{k=0}^{I_1}\varrho _kz=
{\sum}_k\zeta _k\Lambda _+^{(-a+1)}\varrho _kz+{\sum}_k(1-\zeta
_k)\Lambda _+^{(-a+1)}\varrho _kz\\
&={\sum}_ku_k+r,\text{ with } u_k=\zeta _k\Lambda _+^{(-a+1)}\varrho
_kz, \;r={\sum}_k(1-\zeta
_k)\Lambda _+^{(-a+1)}\varrho _kz.
\endaligned
%\tag2.11b
$$
Since $(1-\zeta  _k)\varrho _k=0$, $(1-\zeta  _k)\Lambda _+^{(-a+1)}\varrho
_k$ is a $\psi $do of order $-\infty $, so it maps $z$ into $C^\infty
({\Bbb R}^n)$; moreover, its symbol in local coordinates is holomorphic
for $\operatorname{Im}\xi _n<0$, so it preserves support in $\comega$.
Hence  $r$ is in $\dot C^\infty (\comega)$,
contained in $\dot H^t_p(\comega)\subset H_p^{(a-1)(t)}(\comega)$ for all
$t$. The integral of $Pf\bar g-f \overline {P^*g}$ over $\Omega $ is
zero when $f$ or $g\in \dot C^\infty (\comega )$, so the contributions
from $r$ are zero. Henceforth we  focus on
the  sum
$$
u'={\sum}_ku_k,
$$
where $u_k$ is supported in $\operatorname{supp}\zeta _k\subset U_i$ and belongs to
$H_p^{(a-1)(s)}(\comega)$. A given $v\in H^{(a-1)(s)}(\comega)$ is
similarly decomposed. (A similar localization step should have been
included in the proof of \cite{G16} Th.\ 4.5.)

Now setting
$$
%u_k=\varrho _ku,\quad v_l=\varrho _lv,\quad 
P_{kl}=\psi _lP\psi _k,\quad P_{kl}^*=\psi _kP^*\psi _l,\tag4.27
$$
we can write, since $\psi _k$ is 1 on $\supp u_k$ and $\supp v_k$, 
$$
\ang{r^+Pu',v'}_{\Omega }-\ang{u',r^+{P^*v'}}_{\Omega }=\sum_{k,l\le
J_0}(\ang{r^+P_{kl}u_k,v_l}_{U_i\cap\Omega }
-\ang{u_k,r^+{P_{kl}^*v}}_{U_i\cap\Omega }).\tag4.28
$$
Recall the notation $U_i\cap\Omega =U_i^+$, $V_i\cap\rnp =V_i^+$.
For each pair $(k,l)$ we treat the term by use of the coordinate map for $U_i$,
$i=i(k,l)$.
Denote by $\underline P_{kl}$ the operator on $V_i\subset {\Bbb R}^n$ that $P_{kl}$
carries over to; its kernel is compactly supported in $V_i\times V_i$. In
detail,
$$ 
\underline P_{kl}=\underline \psi _l^{(i)}\underline
P^{(i)}\underline \psi _k^{(i)},\tag4.29$$
 cf.\ (4.21). The
parity  property of the symbol, hence the $a$-transmission property
at any boundary,
is preserved
under the coordinate transformation.

In the $i$'th term, the
two sides are compactly supported in $U_i$. We shall need to carry
this over to a (sesquilinear) distribution duality over $V_i$ using the coordinate
change $\kappa _i$, and therefore recall some general rules:
$$\aligned
\ang {f,g}_{U_i}&=
\ang{\underline f,J\underline g }_{V_i}\\
\ang{Tf,g}_{U_i}&=\ang{f,T^*g}_{U_i}=\ang{\underline f,J\underline{T^*g}}_{V_i}\\
\ang{Tf,g}_{U_i}&=\ang{\underline T\underline f,J\underline
g}_{V_i}=\ang{\underline f,(\underline T)^{(*)}(J\underline
g)}_{V_i}\\
&=\ang{f,J^{-1}[(\underline T)^{(*)}(J\underline g)]\circ \kappa _i}_{U_i}.
\endaligned$$
Here the underlined objects are the elements carried over to $V_i$, and
$J$ is the Jacobian of the mapping $\kappa _i^{-1}$ (the absolute value of
its functional determinant). $J$ is a smooth positive function; for
simplicity of notation we
leave out underlines and the marking of $i$-dependence there.
The star indicates the adjoint with respect to integration over $U_i$, and the
star in parentheses indicates the adjoint with respect to integration
over $V_i$.
We see that the two concepts of adjoints are related by 
$$
J\underline{T^*g}=(\underline T)^{(*)}(J\underline g).
$$

We have for the $(k,l)$'th term:
$$
\ang{r^+P_{kl}u_k,v_l}_{U_i^+}=\ang{r^+\underline P_{kl}\underline
u_k,J\underline v_l}_{V_i^+},\quad
\ang{u_k,r^+P_{kl}^*v_l}_{U_i^+}=\ang{\underline
u_k,r^+J\underline{(P_{kl}^*)}\underline v_l}_{V_i^+},
$$
so with $\underline w_l=J\underline v_l$,
$$
\ang{r^+P_{kl}u_k,v_l}_{U_i^+}-\ang{u_k,r^+P_{kl}^*v_l}_{U_i^+}=\ang{r^+\underline P_{kl}\underline
u_k,\underline w_l}_{V_i^+}-\ang{\underline
u_k,r^+(\underline P_{kl})^{(*)}\underline w_l}_{V_i^+},\tag4.30
$$
where $(\underline P_{kl})^{(*)}$ is the adjoint in the $V_i^+$-setting,
$$
(\underline P_{kl})^{(*)}=J\underline{(P_{kl}^*)}J^{-1}.
$$
The underlined operators are defined on $V_i$, where $\kappa
_i(U_i\cap \comega )\subset\crnp $ and the relevant part 
of the boundary of $V_i\cap \crnp$ is a subset of ${\Bbb R}^{n-1}$ (all objects are
supported away from the other parts of the boundary). 

We can consider $\underline P^{(i)}$ as extended to an operator of the
same type on all of ${\Bbb R}^n$ with global estimates (and keep the name  $\underline P^{(i)}$), so
that (4.29) holds with the globally defined operator; this is followed
up for the adjoint.

Now we are in a situation to apply Theorem 4.1, which gives
$$\aligned
&\ang{r^+\underline P_{kl}\underline
u_k,\underline w_l}_{V_i^+}-\ang{\underline
u_k,r^+(\underline P_{kl})^{(*)}\underline w_l}_{V_i^+}=\ang{r^+\underline P_{kl}\underline
u_k,\underline w_l}_{\rnp}-\ang{\underline
u_k,r^+(\underline P_{kl})^{(*)}\underline w_l}_{\rnp}\\
&\quad =(\underline s_{kl,0}\gamma _1^{a-1}\underline u_k,\gamma
_0^{a-1}\underline w_l)-(\underline s_{kl,0}\gamma _0^{a-1}\underline
u_k,\gamma _1^{a-1}\underline w_l)+(\underline B_{kl}\gamma
_0^{a-1}\underline u_k,\gamma _0^{a-1}\underline w_{l}),
\endaligned\tag4.31$$
the last line consists of ${L_2({\Bbb R}^{n-1})}$-scalar
products. Here $\underline s_{kl,0}=\underline \psi _l \underline
p_0(x',0,0,1)\underline \psi _k$.

The dualities over $V_i^+$ carry over to dualities over $U_i^+$ by (4.30).
For the scalar products over ${\Bbb R}^{n-1}$ (supported in $V_i'$),
we note that the coordinate transform preserves the definitions of
$\gamma _0^{a-1}$ and $\gamma _1^{a-1}$ (since $t$ corresponds exactly
to $x_n$ for $t<r'$),
$$
[\gamma _0^{a-1}\underline u_k]\circ \kappa _i'=\gamma _0^{a-1}u_k,
\quad [\gamma _1^{a-1}\underline u_k]\circ \kappa _i'=\gamma _1^{a-1}u_k.$$
Moreover,
$$\aligned
\gamma _0^{a-1}\underline w_l&=\Gamma (a)\gamma _0(J\underline
v_l/t^{a-1})=J_0\gamma _0^{a-1}\underline v_l,\\
\gamma _1^{a-1}\underline w_l&=\Gamma (a+1)\gamma _0(\partial_t(J\underline
v_l/t^{a-1}))=J_0\gamma _1^{a-1}\underline v_l+aJ_1\gamma
_0^{a-1}\underline v_l,
\endaligned$$ 
 where $J_0=\gamma _0(J)$, $J_1=\gamma _0(\partial_tJ)$. Here $J_0$
 defines the area element $d\sigma $ in integrations over the boundary:
$$
\int_{U'_i}f(x')\, d\sigma =\int_{V'_i}f((\kappa '_i)^{-1}(y'))J_0\,dy',
$$
and $J_1$ gives rise to an extra term along with  $\gamma
_1^{a-1} v_l$.
Then (4.31) carries over to the
formula on $U_i^+$:
$$\aligned
&\ang{r^+P_{kl}u_k,v_l}_{U_i^+}-\ang{u_k,r^+P_{kl}^*v_l}_{U_i^+}=( s_{kl,0}\gamma _1^{a-1}u_k,\gamma
_0^{a-1}v_l)_{U'_i}-( s_{kl,0}\gamma _0^{a-1}u_k,\gamma
_1^{a-1}v_l)_{U'_i}\\
&\quad+((B_{kl}-as_{kl,0}J_0^{-1}J_1)\gamma
_0^{a-1}u_k,\gamma _0^{a-1}v_l)_{U'_i},
\endaligned\tag4.32$$
where $s_{kl,0}=\psi _lp_0(x,\nu (x))\psi _k$.
Here $B_{kl}$ depends only on the principal and subprincipal symbols of
$P_{kl}$, cf.\ Remark 4.2.

Finally, summing over $k$
and $l$ and using that 
$\gamma _0^{a-1}u=\gamma _0^{a-1}u'=\sum_k\gamma _0^{a-1}u_k$,
%$\psi _k\varrho _k=\varrho _k$, 
we find (4.25)
with 
$$
B\gamma _0^{a-1}u =\sum_{k,l}(B_{kl}-as_{kl,0}J_0^{-1}J_1)\gamma
_0^{a-1}u_k.
%B \gamma _0^{a-1}u =\sum_{k,l}\varrho _l[(\underline
%B_{kl}-a\underline s_{kl,0}J_0^{-1}J_1)\underline{\varrho
%_k\gamma ^{a-1}_0u}]\circ \kappa _{i(k,l)}'.
\tag4.33 
$$

The last statements are seen as in Theorem 3.3.
\qed

\enddemo

\proclaim{Corollary 4.5} Let $P$ and $\Omega $ be as in Theorem {\rm
4.3}.

$1^\circ$ If $u\in H^{(a-1)(s)}(\comega)$, $v\in
H^{a(s)}(\comega)$,  $s>a+\frac12$, then
$$
\aligned
\ang{r^+Pu,v}&_{\ol H^{-a+\frac12+\varepsilon }(\Omega ), \dot
H^{a-\frac12-\varepsilon }(\comega)}-\ang{u,r^+P^*v}_{\dot H^{a-\frac12-\varepsilon }(\comega),\ol
H^{-a+\frac12+\varepsilon }(\Omega )}\\
&=-( s_0\gamma _0^{a-1}u,\gamma _0^av)_{L_2(\partial\Omega )}\\
&=-\Gamma (a)\Gamma (a+1)\int_{\partial\Omega }s_0\gamma _0(\tfrac{u}{t^{a-1}})\gamma
_0(\tfrac{\bar v}{t^{a}})
\, d\sigma .
\endaligned\tag4.34
$$

$2^\circ$ It follows that if  $u,v\in
H^{a(s)}(\comega)$,  $s>a+\frac12$, then for each $j=1,\dots,n$,
$$
\aligned
\ang{r^+Pu,\partial_jv}&_{\ol H^{-a+\frac12+\varepsilon }(\Omega ), \dot
H^{a-\frac12-\varepsilon }(\comega)}
+\ang{\partial_ju,r^+P^*v}_{\dot H^{a-\frac12-\varepsilon }(\comega),\ol
H^{-a+\frac12+\varepsilon }(\Omega )}\\
&=(\nu _j s_0\gamma _0^au,\gamma _0^av)_{L_2(\partial\Omega
)}+(r^+[P,\partial_j]u,v)_{L_2(\Omega )}\\
&=\Gamma (a+1)^2\int_{\partial\Omega }\nu _j s_0\gamma _0(\tfrac
u{t^a})\,\gamma _{0}(\tfrac {\bar v}{t^a})\, d\sigma +(r^+[P,\partial_j]u,v)_{L_2(\Omega )}.
\endaligned \tag4.35
$$

The dualities in the left-hand sides can be written as integrals when
$s\ge 2a$.
\endproclaim 

\demo{Proof} $1^\circ$ follows simply by application of Theorem 4.4 with $\gamma
_0^{a-1}v=0$; then $\gamma _1^{a-1}v=\gamma _0^av$, as noted
earlier.

$2^\circ$ is deduced from this as follows: 
First we observe that 
$$
\ang{r^+Pu,\partial_jv}_{\ol H^{-a+\frac12+\varepsilon }, \dot
H^{a-\frac12-\varepsilon }}=-\ang{\partial_jr^+Pu,v}_{\ol H^{-a-\frac12+\varepsilon }, \dot
H^{a+\frac12-\varepsilon }}
$$
by integration by parts, using that $\gamma _0v=0$ ($v\in
H^{a(a+\frac12+\varepsilon )}(\comega)\subset\dot
H^{a+\frac12-\varepsilon }(\comega)$). Next, we introduce the commutator 
$[P,\partial_j]=P\partial_j-\partial_jP
$. Now $\partial_ju\in
H^{(a-1)(s)}$, and (4.34) applies:
$$
\aligned
&\ang{r^+Pu,\partial_jv}_{\ol H^{-a+\frac12+\varepsilon }, \dot
H^{a-\frac12-\varepsilon }}
+\ang{\partial_ju,r^+P^*v}_{\dot H^{a-\frac12-\varepsilon },\ol
H^{-a+\frac12+\varepsilon }}\\
&=-\ang{r^+P\partial_ju,v}_{\ol H^{-a-\frac12+\varepsilon }, \dot
H^{a+\frac12-\varepsilon }}
+\ang{\partial_ju,r^+P^*v}_{\dot H^{a-\frac12-\varepsilon },\ol
H^{-a+\frac12+\varepsilon }}+(r^+[P,\partial_j]u,v)\\
&=( s_0\gamma _0^{a-1}(\partial_ju),\gamma _0^av)+(r^+[P,\partial_j]u,v).
\endaligned\tag4.36
$$
Write $u$ near $\partial\Omega $ as $u=t^aw(y'+\nu (y')t)$, where
$y'\in\partial\Omega $ and $w$ is constant in $t$ for $0\le t<r'$, and
use that $\partial_j=\nu _j(y')\partial_t+T$, where $T$ is tangential
(acts along $\partial\Omega $) there, to see that
$$
\partial_ju=\nu _j(y')\partial_t(t^aw)+T(t^aw)=\nu _jat^{a-1}w+t^aTw,
$$
 hence
$$
\aligned
\gamma _0^{a-1}(\partial_ju)&=\Gamma (a)\gamma _0(\partial
_ju/t^{a-1}))
=\Gamma (a)\gamma _0(\nu _jaw+tTw)\\
&=\Gamma (a+1)\nu _j\gamma _0w=\Gamma (a+1)\nu _j\gamma _0(u/t^a)=\nu
_j\gamma _0^au.
\endaligned
$$
Insertion in (4.36) leads to (4.35). \qed

\enddemo

A version of $1^\circ$ in the corollary was shown by Abatangelo
\cite{A15}, see (9), in the
case $P=(-\Delta )^a$, $0<a<1$. Since $d^{a-1}$ blows up for $d\to 0$
in this case, the functions $u$ entering in the formula are called
``large solutions'' in \cite{A15}.

$2^\circ$ was shown by Ros-Oton and Serra for $(-\Delta )^a$, $0<a<1$, in
\cite{RS14b} (under different smoothness hypotheses), extended to
higher-order fractional Laplacians in \cite{RS15}, and generalized to
larger classes of translation-invariant operators
in a joint work with Valdinoci \cite{RSV17}. We extended it to
($x$-dependent) elliptic $\psi $do's  in \cite{G16}. It implies
Pohozaev formulas that are used to show uniqueness results for
nonlinear problems. 

Note that no assumption on ellipticity of $P$ is made in the theorem
and corollary. Actually, it is not surprising that ellipticity is not needed, since the formulas
are linear in $P$: For a given $P$ one can add $c(1-\Delta )^a$ with a
sufficiently large constant $c$ to make the sum $P_c$ strongly
elliptic; a result for $P_c$ will then lead to a result for $P$ by
subtraction of the formula for $c(1-\Delta )^a$. However, ellipticity
was an important ingredient in earlier proofs of the corollary, that
we can now do away with.

The corollary only involves first boundary values, and the boundary
contribution is purely local. Theorem 4.4, however, lets the Neumann
value enter in a nontrivial way along with a nontrivial Dirichlet
value; this seems to be entirely new even for $(-\Delta )^a$.
\medskip

%\endexample

\subhead 5. A parametrix of the Dirichlet problem 
\endsubhead

Green's formula shows that the Dirichlet and Neumann trace operators
$\gamma _0^{a-1}$ and $\gamma _1^{a-1}$ play a fundamental role in the
discussion of boundary value problems for $P$. Much is known about the
homogeneous Dirichlet problem (1.5), whereas problems for $P$ with
nonzero Dirichlet trace have been less studied. We showed the Fredholm
solvability in \cite{G15, G14} for large scales of spaces over a
smooth bounded set $\Omega $. The structure of the solution operator
will now be further clarified, in a study of the operator $K_D$
solving (1.9).

In addition to the assumptions listed in the beginning of Section 4,
we assume that $P$ is elliptic, 
the principal symbol avoiding a ray. 
This holds in particular if $P$ is strongly elliptic, i.e.
$$ 
\operatorname{Re}p_0(x,\xi )\ge c|\xi |^{2a} \text{ for }|\xi |\ge
1, \text{ with }c>0.\tag5.1
$$
%it is said to be uniformly
% strongly elliptic on ${\Bbb R}^n$ when $c(x)$ can be taken constant.
%More generally, factorization index $a$ holds when $P$ is elliptic avoiding a
% ray, cf.\ \cite{G16} Sect.\ 2.4. 
 
The following mapping was defined in \cite{G15} for
$s>a-\frac12$,% $\Omega $ bounded,
$$
\pmatrix r^+P\\
\gamma _0^{a-1}\endpmatrix\colon H^{(a-1)(s)}(\comega)\to \ol
H^{s-2a}(\Omega )\times H^{s-a+\frac12}(\partial\Omega );\tag5.2
$$
by Th.\ 6.1 there it is Fredholm when $\Omega $ is bounded.
Extensions of the mapping property to $H^s_p$-spaces and other scales of spaces are given
in \cite{G14, G15}. To exhibit a bijective case, we show:

\proclaim{ Lemma 5.1} 
Let $P$  be a classical $\psi $do on ${\Bbb
R}^n$ of order $2a>0$, strongly elliptic  with even symbol on ${\Bbb
R}^n$, and let
$\Omega $ be a smooth bounded subset. 
If
$P$ has a positive
lower bound on $C_0^\infty (\Omega )$:
$$
\operatorname{Re}(Pu,u)_{L_2(\Omega )}\ge c_0\|u\|_{L_2(\Omega )}^2\text{ for }u\in
C_0^\infty (\Omega ),\tag 5.3
$$
with $c_0>0$, then the mapping {\rm (5.2)} is bijective. 
 The solution operator, denoted $\pmatrix R_D &
K_D\endpmatrix$, maps as follows:
$$
\pmatrix r^+P\\ \gamma _0^{a-1}\endpmatrix^{-1}=\pmatrix R_D &
K_D\endpmatrix\colon  \ol
H^{s-2a}(\Omega )\times H^{s-a+\frac12}(\partial\Omega )\to
H^{(a-1)(s)}(\comega),\; s>a-\tfrac12.\tag 5.4
$$
\endproclaim

\demo{Proof}
As mentioned in Section 1, the Dirichlet realization
$P_D$ defined by the variational construction
(the Lax-Milgram lemma) acts
like $r^+P$ with domain $$
D(P_D)=\{u\in \dot H^a(\comega)\mid r^+Pu\in
L_2(\Omega )\}.\tag5.5
$$
It was found in \cite{G15} that $D(P_D)=H^{a(2a)}(\comega)$. The
adjoint is the analogous operator for $P^*$, which also satisfies
(5.3). This inequality assures that both $P_D$ and $P_D^*$ are
injective, hence $P_D$ is bijective from $H^{a(2a)}(\comega)$ to
$L_2(\Omega )$.
Now
$H^{a(2a)}(\comega)=\{u\in H^{(a-1)(2a)}(\comega)\mid \gamma
_0^{a-1}u=0 \}$, and 
$$\gamma _0^{a-1}\colon H^{(a-1)(2a)}(\comega)/ H^{a(2a)}(\comega)\simto
H^{a+\frac12}(\partial\Omega )\tag5.6
$$
(cf.\ \cite{G15}, Th.\ 5.1),
 so the bijectiveness in (5.2) for $s=2a$ follows by supplying the mapping 
$P_D$ with (5.6).
 The bijectiveness holds also for
other $s$ in view of the invariance of kernels and cokernels
(\cite{G14}, Th\. 5.5).\qed

\enddemo

An example of an operator satisfying the hypotheses of Lemma 5.1 is
$(1-\Delta )^a$, since 
$$
((1-\Delta )^au,u)=(2\pi )^{-n}\int_{{\Bbb R}^n}(1+|\xi |^2)^{2a}|\hat u(\xi )|^2
\,d\xi 
%\ge (2\pi )^{-n}\int_{{\Bbb R}^n}|\hat u(\xi )|^2
%\,d\xi 
\ge\|u\|^2_{L_2({\Bbb R}^n)} \text{ for all }u\in \Cal S({\Bbb R}^n) .$$
In view of (3.12)ff.\ and (3.14), the operator $K_D$ for
$P=(1-\Delta )^a$ on $\rnp$ equals $K^{a-1}_0$.

In \cite{G15} we obtained in Th.\ 4.4 that 
$$
R_D=\Lambda _+^{(-a)}e^+\widetilde{Q_+}r^+\Lambda ^{(-a)}_{-}e^+\colon  \ol
H^{s-2a}(\Omega )\to H^{(a)(s)}(\comega)\tag5.7
$$
is a parametrix of the Dirichlet
problem with zero boundary condition; here $Q=$ \linebreak$\Lambda
_-^{(-a)}P\Lambda _+^{(-a)}$,
$Q_+=r^+Qe^+\colon \ol H^t(\Omega)\to \ol H^t(\Omega)$ (for all $t>-\frac12$), and $\widetilde{Q_+}$ is a parametrix of $Q_+$.
Moreover, we showed in Th.\ 6.5 there how, in the case $\Omega =\rnp$,
$R_D$ could be supplied with a Poisson-like operator $$
K_D\colon H^{s-a+\frac12}({\Bbb R}^{n-1})\to
H^{(a-1)(s)}(\crnp)$$
constructed from $R_D$, to give a full parametrix. The operator $K_D$ was shown to be of the form
$$
K_D=\Xi _+^{1-a}e^+K'=\Lambda  _+^{1-a}e^+K'',
$$
with $K'$ and $K''$ being Poisson operators of order 0 belonging to the Boutet de
Monvel calculus. 
%They were described there in terms of $r^+Qe^+$ and
%its parametrix, where $Q=\Xi _-^{-a}P\Xi _+^{-a}$.

In \cite{G16} Th.\ 2.7, assuming that $P$ is elliptic avoiding a ray,
with even symbol, we worked out an approximate factorization of $P$ in "minus" and
"plus" operators (preserving support in $\crnm$ resp.\ $\crnp$)   
in the case $\Omega =\rnp$, 
$$
P\sim P^-P^+,%\text{ with 
%parametrices }\widetilde P, \widetilde Q^-, \widetilde Q^+,\text{ such
%that } \widetilde Q\sim \widetilde Q^+ \widetilde Q^-,
\tag5.8$$ 
and we indicated  in Rem.\
2.9 there  a formula for $K_D$ based on the factorization, namely essentially 
$$
K_D \varphi  \sim r^+\widetilde P^+( \varphi  (x')\otimes \delta (x_n)),\tag5.9
$$ 
where $\widetilde P^+$ is a parametrix of $P^+$. 
We shall now go into details with this construction, and thereby also obtain a more
informative
formula for the full parametrix $\pmatrix R_D&K_D\endpmatrix$.

As accounted for in the Appendix,
we have the following
%Returning to $P$, this leads to 
product decompositions of $P$ and its
parametrix $\widetilde P$:
$$
\aligned
P&=\Xi _-^aQ\Xi _+^a=\Lambda _-^aQ_1\Lambda _+^a\sim P^-P^+,\text{
where}\\
Q&=\Xi _-^{-a}P\Xi
_+^{-a}\sim Q^-Q^+, \quad Q_1=\Lambda  _-^{-a}P\Lambda 
_+^{-a}\sim Q_1^-Q_1^+,\text{ and}\\
\quad P^-&=\Xi _-^aQ ^-=\Lambda _-^aQ_1^-,
\quad
P^+
=Q ^+\Xi _+^a=Q_1^+\Lambda _+^a;\text{ moreover,}\\
\widetilde P&\sim \widetilde P^+\widetilde P^-,\quad \widetilde
P^+=\Xi _+^{-a}\widetilde Q ^+=\Lambda _+^{-a}\widetilde Q_1^+
,\quad
\widetilde P^-
=\widetilde Q ^-\Xi _-^{-a}=\widetilde Q_1^-\Lambda _-^{-a}.
\endaligned\tag5.10
$$
In view of the support-preserving properties, when $\supp
u\subset \crnp$, $P^+u\in \dot H^{-\frac12+\varepsilon }(\crnp)$,
$$
%\aligned
r^+Pu\sim r^+P^-P^+u=r^+P^-e^+r^+P^+u,
%r^+Q_1u&\sim r^+Q_1^-Q_1^+u=r^+Q_1^-e^+r^+P^+u.
%\endaligned
\tag5.11
$$
with similar rules for $Q$ and $Q_1$.

Let us analyse how $K_D$ should look, by an argumentation that reduces
the problem as far as possible to the Boutet de Monvel calculus, in the same
spirit as the proof of Th.\ 4.4 in \cite{G15}. (One advantage of that
calculus is that there are good rules for smoothing operators, which
can be difficult to obtain in a mixture of fractional-order $\psi $do's
and the generalized $\psi $do's that in some situations only give good
tangential results.)

We want to solve:
$$
r^+Pu=0,\quad \gamma _0^{a-1}u=\varphi ,
\tag5.12
$$
for a given $\varphi \in H^{s-a+\frac12}({\Bbb R}^{n-1})$ with $s>a-\tfrac12$, searching for
$u$ in $H^{(a-1)(s)}(\crnp)$. 
Let $\tilde u\in \ol H^{s-a+1}(\rnp)$ be
the function for which
$$
\tilde u=r^+\Xi _+^{1-a}u;
 \text{ then }u=\Xi _+^{a-1}e^+\tilde u\text{ and } e^+\tilde u=\Xi _+^{1-a}u,\tag5.13
$$
by definition (cf.\ \cite{G15}, Sect.\ 1). For example when $s=a$,
$\varphi \in H^{\frac12}({\Bbb R}^{n-1})$ and  $\tilde u\in \ol
H^1(\rnp)$. 

By use of (5.10), (5.11), the equations in (5.12) can be written
$$
r^+\Xi _-^ae^+r^+Q \Xi _+^1\Xi _+^{a-1}u=0,\quad \gamma _0\Xi _+^{a-1}u=\varphi .
$$
Using (5.13), this may, since $r^+\Xi _-^ae^+$ is bijective,
equivalently be written
$$
r^+Q \Xi _+^1e^+\tilde u=0,\quad \gamma _0\tilde u=\varphi ,\tag 5.14
$$
where the boundary value is as always taken from $\rnp$. 

Now (5.14) is concerned with integer-order operators, where we have the rules of the Boutet de Monvel
calculus, extended to allow the generalized $\psi $do's $Q^+$ and $\Xi _+^1$.
Here we note that %$Q \sim Q^-Q^+$ where $Q^\pm$ have parametrices $\widetilde Q^\pm$,and 
$r^+Q \Xi _+^1e^+\tilde u\sim r^+Q ^-e^+r^+Q ^+\Xi
_+^1e^+\tilde u$, where $ r^+Q ^-e^+$ has the parametrix $
r^+\widetilde Q ^-e^+$. This allows transforming
(5.14) to 
$$
r^+Q ^+\Xi _+^1e^+\tilde u\sim 0,\quad \gamma _0\tilde u=\varphi ,\tag 5.15
$$
in a parametrix sense. 
Here, since   $Q ^+\Xi _+^1$ has the parametrix
$\Xi _+^{-1}\widetilde Q ^+$, the problem (5.15) can be expected to have the solution operator
$$
K_{Q^+ }\varphi \equiv r^+\Xi _+^{-1}\widetilde Q ^+(\varphi (x')\otimes \delta (x_n)),\tag5.16
$$
in a parametrix sense. This is the Poisson operator with symbol
$k_{Q^+}\in S^{-1}(\Cal H^+)$,
$$
k_{Q^+ }(x',\xi )=\chi _+^{-1}(\xi )\#{\widetilde q^{+\prime}}(x',0,\xi ),\tag5.17
$$
by the rule (A.13), the prime indicating that we have used the 
 $(x',y_n)$-form of the symbol of $\widetilde Q^+$, cf.\ (A.11). It is known from the calculus that $K_{Q^+}$ maps
$$
K_{Q^+}\colon H^{t-\frac12}({\Bbb R}^{n-1})\to \ol H^{t}(\rnp)\text{ for }t\in{\Bbb R}.\tag5.18
$$
In case $Q=I$, $K_{Q^+ }$ is the standard
 Poisson operator $K_0=\operatorname{OPK}(\chi _+^{-1})$ which
 satisfies $\gamma _0K_0=I$.

With this solution operator to (5.15), we go on to define the
solution operator for (5.12) by
$$
K_D=\Xi _+^{1-a}e^+K_{Q^+}\colon H^{s-a+\frac12}({\Bbb R}^{n-1})\to \Xi
_+^{1-a}e^+ \ol H^{s-a+1}(\rnp)=H^{(a-1)(s)}(\crnp).\tag5.19
$$
It is accounted for how this $K_D$ has the desired parametrix property in the proof of the following theorem.

\proclaim{Theorem 5.2} Let $P$ be a classical globally estimated $\psi
$do  on ${\Bbb R}^n$ of order
$2a>0$, elliptic avoiding a ray,  with even
symbol. 

$1^\circ$ Define $K_D$ on
$H^{s-a+\frac12}({\Bbb R}^{n-1})$, $s>a-\frac12$, by {\rm (5.19)}
with {\rm (5.16)}.
It  solves the Dirichlet problem
$$
r^+Pu=0\text{ on }\rnp,\quad  \gamma _0^{a-1}u=\varphi \text{ at
}x_n=0,\tag 5.20
$$
in a parametrix sense, namely
$$
 \gamma _0^{a-1}K_D=I, \quad r^+PK_D=\Cal S\colon  H^{s-a+\frac12}({\Bbb R}^{n-1})\to
C^\infty (\crnp).\tag5.21
$$

$2^\circ$ Define instead $K_D$ by
$$
K_D=\Lambda  _+^{1-a}e^+K_{Q_1^+}.\text{ where }
K_{Q_1^+}\varphi =r^+\Lambda  _+^{-1}\widetilde Q_1^+(\varphi (x')\otimes
\delta (x_n)).
\tag5.22
$$
Then we again have that $K_D$ solves {\rm (5.20)} in a parametrix
sense, namely {\rm (5.21)} holds (with a possibly different $\Cal S$).
%it is defined in {\rm (5.11)} below.

\endproclaim

\demo{Proof}
$1^\circ$.
 We begin by checking that $K_{Q^+}$ defined in (5.16)  solves the
 problem (5.15). First,
$$
r^+Q ^+\Xi _+^{1}e^+K_{Q^+}\varphi =r^+Q ^+\Xi _+^{1}e^+r^+\Xi
_+^{-1}\widetilde Q^+(\varphi \otimes \delta )= r^+Q ^+\widetilde Q ^+(\varphi
\otimes\delta )=\Cal S_1\varphi ,
$$
where $\Cal S_1$ is a Poisson operator of order $-\infty $, hence maps
$H^{s-a+\frac12}({\Bbb R}^{n-1})$ into $C^\infty (\crnp)$ (more
precisely, into the subset $\bigcap_t\ol H^t(\rnp)$). This uses that
$(I-e^+r^+)v=0$ when $v=\Xi
_+^{-1}\widetilde Q^+(\varphi \otimes \delta )$, since it is an
$L_2$-function supported in $\crnp$; then since $Q ^+\widetilde
Q ^+=I+\Cal R$, where $\Cal R$ has symbol  $r(x',\xi )\in S^{-\infty }(\Cal
H^+)$, $\Cal S_1$ equals $\operatorname{OPK}(r)$, a smoothing Poisson operator.

Next, since $\widetilde Q^+$ has a symbol of the form $1+f^+$, where
$f^+=h^+\widetilde q^+\in S^0(\Cal H^+)$,
$$
\gamma _0K_{Q^+}\varphi 
%\Xi _+^{-1}\widetilde Q^+_1(\varphi \otimes \delta )
=\gamma _0\Xi
_+^{-1}(I+F^+)(\varphi \otimes \delta )=\gamma _0K_0\varphi
 +\gamma _0\Xi _+^{-1}F^+(\varphi \otimes\delta )=\varphi ,
$$
where the term with $\Xi _+^{-1}F^+$ gives zero since the symbol is in
$\Cal H^+$ w.r.t.\ $\xi _n$ and is $O(\xi _n^{-2})$, so that its
integral in $\xi _n$ vanishes.
Altogether,
$$
r^+Q ^+\Xi _+^{1}e^+K_{Q^+}=\Cal S_1,\quad \gamma _0K_{Q^+}=I.\tag 5.23
$$

Now define $K_D$ by (5.19). Then
$$
r^+PK_D\varphi =r^+\Xi _-^ae^+r^+Q \Xi _+^aK_D\varphi =
r^+\Xi _-^ae^+r^+(Q ^-Q ^++\Cal R_1)\Xi _+^1e^+K_{Q^+}\varphi,\tag5.24
$$
where $\Cal R_1$ has symbol in $S^{-\infty} (\Cal H_{-1})$.
Here 
$$
r^+Q ^-Q ^+\Xi _+^1e^+K_{Q^+}\varphi=r^+Q ^-e^+r^+Q ^+\Xi
_+^1e^+K_{Q^+}\varphi=r^+Q ^-e^+\Cal S_1\varphi=\Cal S_2\varphi ,
$$
in view of (5.23), where $\Cal S_2$ maps into $\bigcap_t\ol
H^t(\rnp)$, since $r^+Q ^-e^+$ preserves this space (by \cite{G16}
Th.\ 2.4 combined with support-preserving properties). Composition with 
$r^+\Xi _-^ae^+$
 to the left gives another smoothing operator $\Cal
S_3$ in view of the isomorphism properties of $r^+\Xi _-^ae^+$.
The contribution from $\Cal R_1$ in (5.24) equals
$$
r^+\Xi _-^ae^+r^+\Cal R_1\Xi _+^1e^+K_{Q^+}\varphi =\Cal S_4\varphi ,
$$
where $\Cal S_4$ is smoothing, since $r^+\Cal R_1\Xi _+^1e^+K_{Q^+}$
is a Poisson operator with symbol in $S^{-\infty}(\Cal H^+) $.

Concerning the boundary condition, we note that
$$
\gamma _0^{a-1}K_D\varphi =\gamma _0\Xi _+^{a-1}
K_D\varphi =
\gamma _0\Xi _+^{a-1}\Xi _+^{1-a}e^+K_{Q^+}\varphi=\gamma _0K_{Q^+}\varphi =\varphi , 
$$
 in view of (5.23). Taking $\Cal S=\Cal S_3+\Cal S_4$, we have obtained
 $1^\circ$.

$2^\circ$. We know from  \cite{G15} that
$H^{(a-1)(s)}(\crnp)$ can equivalently be defined as \linebreak $\Lambda
_+^{1-a}e^+\ol H^{s-a+1}(\rnp)$. Moreover, by  (A.16),
$$
\gamma _0^{a-1}u=\gamma
_0(\Xi ^{a-1}_+u)= \gamma _0(\Lambda  ^{a-1}_+u).\tag5.25
$$
The Poisson operator $K_0'=\operatorname{OPK}(\lambda _+^{-1})$
satisfies $\gamma _0K'_0=I$, and $\Lambda _+^{1-a}e^+K'_0$ is a
right-inverse of $\gamma _0^{a-1}$. 

The  proof under $1^\circ$ goes over verbatim to a proof of $2^\circ$ when the $\Xi
_\pm^t$-family is replaced by $\Lambda _\pm^{t}$, $Q $ is replaced by
$Q_1$, and $\tilde u$ is replaced by $\tilde u'=r^+\Lambda _+^{a-1}u$.
\qed

\enddemo

\proclaim{Corollary 5.3}  Together with $R_D$ recalled above, $K_D$  enters in a full parametrix
 $$
\pmatrix  R_D&K_D\endpmatrix\colon \ol H^{s-2a}(\rnp)\times
H^{s-a+\frac12}({\Bbb R}^{n-1})\to H^{(a-1)(s)}(\crnp),\quad s>a-\tfrac12,\tag5.26
$$
satifying$$
\pmatrix r^+P\\ \gamma _0^{a-1}\endpmatrix \pmatrix R_D & K_D\endpmatrix
=\pmatrix I& 0\\ 0&I\endpmatrix+\pmatrix \Cal S_1&\Cal S\\0&0\endpmatrix,\tag5.27
$$ 
where $\Cal S_1$ maps $\ol H^{s-2a}(\rnp)$ into $C^\infty (\crnp)$ and
$\Cal S$ is as in {\rm (5.21)}.
\endproclaim

\demo{Proof} From the way in which $R_D$ was defined in \cite{G15}, we  have
that $\gamma _0^{a-1}R_D=0$ and $r^+PR_D=I+\Cal S_1$ where $\Cal S_1$ is a
smoothing operator; (5.27) follows by combining this with the
above theorem.\qed
\enddemo

By use of local coordinates (cf.\ Remark 4.3), the above construction of $K_D$ in the case
$\rnp$ can be applied to construct the Poisson-like operator in the
general case of a bounded smooth open set  $\Omega $.
%Here we shall use the coordinate charts described in Remark 4.3.

\proclaim{Theorem 5.4} Let $P$ be a classical $\psi $do of order
$2a>0$, elliptic avoiding a ray,  with even
symbol, and let $\Omega $ be a smooth bounded subset of ${\Bbb
R}^n$. 
There is a Poisson-like operator $K_D\colon H^{s-a+\frac12}(\partial\Omega )\to
H^{(a-1)(s)}(\comega) $,  $s>a-\frac12$ (see details in {\rm (5.31)}), such that together with
$R_D$ constructed in \cite{G15} {\rm Th.\ 4.4}, $$
\pmatrix
R_D&K_D\endpmatrix\colon \ol H^{s-2a}(\Omega )\times
H^{s-a+\frac12}(\partial\Omega )\to H^{(a-1)(s)}(\comega),\quad s>a-\tfrac12,\tag5.28
$$ is a parametrix of the nonhomogeneous Dirichlet
problem, satisfying $$
\pmatrix r^+P\\ \gamma _0^{a-1}\endpmatrix \pmatrix R_D & K_D\endpmatrix
=\pmatrix I& 0\\ 0&I\endpmatrix+\pmatrix \Cal S_1&\Cal S\\0&0\endpmatrix,\tag5.29
$$ 
where $\Cal S_1\colon \ol H^{s-2a}(\Omega )\to C^\infty (\comega)$ and
$\Cal S\colon  H^{s-a+\frac12}(\partial\Omega )\to
C^\infty (\comega)$.
\endproclaim

\demo{Proof} 
%In order to factorize $P$ in factors of minus- and
%plus-type, 
We shall here use the cover  $U_i$, $i=0,\dots, I_0$, described
in Remark 4.3, with an {\it associated}  partition of unity $\varrho 
_0,\dots,\varrho _{I_0}$ such that each $\varrho _i$ is in $ C_0^\infty (U_i)$
taking values in $[0,1]$, and $\sum_{0\le i\le I_0}\varrho _i(x)=1$ on
a neighborhood of $\comega$. 
One can moreover find functions
 $\zeta _i^0,\zeta _i^1,\zeta _i^2,\dots\in C_0^\infty
(U_i,[0,1])$ such that $\zeta _i^{0}\varrho _i=\varrho _i$ and $\zeta
_i^{k+1}\zeta _i^k=\zeta _i^k$ for all $k$ (in each case, the former function is 1
on the support of the latter).

$P$ considered on functions supported in $U_i$ gives rise to
$\underline P^{(i)}$ defined on functions supported in $V_i$ by
(4.21). The symbol of the operator $\underline P^{(i)}$ can be assumed to be extended to a
$\psi $do symbol on ${\Bbb R}^n$, even and of order $2a$, and elliptic
avoiding a ray there (we use the same notation for the extension).
There are factorizations, as described in (5.10), with notation 
$$
\underline P^{(i)}=\Lambda _-^a\underline Q_1^{(i)}\Lambda _+^a, \quad
\underline Q_1^{(i)} \sim  \underline Q_1^{(i)-}  \underline Q_1^{(i)+},
\text{ etc.}\tag5.30
$$

For $\varphi $ given on $\partial\Omega $, $\varphi
=\sum_{0\le i\le I_0}\varrho _i\varphi $, where $\varrho _i\varphi $
carries over to $\underline{\varrho _i\varphi }$ on $V_i'$. Define
$$
\aligned
 K^{(i)}_{\underline Q_1^+}\psi &=r^+\Lambda _+^{-1}\widetilde {\underline Q_1}^{(i)+}(\underline\psi (x')
\otimes\delta (x_n))\text{ when }\supp\underline\psi \subset V'_i,\\
\underline K^{(i)}_D\underline \psi 
&=\Lambda _+^{1-a}e^+ K^{(i)}_{\underline Q_1^+}\underline\psi ,\\
K_D\varphi &=\sum_i(\underline\zeta _i^1\underline K_D^{(i)}\underline{\varrho _i\varphi })\circ \kappa _i.
\endaligned\tag5.31
$$
(The contribution from $U_0$ is 0.) We shall verify that this
definition of $K_D$ leads to the desired properties.

For the composition with $r^+P$ we have:
$$
\aligned
r^+PK_D\varphi &=\sum_ir^+P(\underline\zeta _i^1\underline K_D^{(i)}\underline{\varrho _i\varphi } )\circ \kappa _i\\
&=\sum_ir^+\zeta _i^2P(\underline\zeta _i^1\underline
K_D^{(i)}\underline{\varrho _i\varphi } )\circ \kappa _i+\Cal S_2\varphi ,
\endaligned
$$
with a smoothing operator $\Cal S_2$, since $1-\zeta _i^2$ and $\zeta _i^1$
have disjoint supports so that $(1-\zeta _i^2)P\zeta _i^1$ is a $\psi
$do of order $-\infty $.

Consider here the $i$'th term carried over to $V_i$, omitting $(i)$
from the notation for simplicity: 
$$
\aligned
&r^+\underline \zeta _i^2\underline P\underline\zeta _i^1 {\underline K}_D\underline{\varrho _i\varphi }
=
r^+\underline \zeta _i^2\Lambda _-^ae^+r^+\underline Q_1\Lambda
_+^a\underline\zeta _i^1\Lambda _+^{1-a}e^+ K_{\underline
Q_1^+}\underline{\varrho _i\varphi }\\
&=
r^+\underline \zeta _i^2\Lambda _-^ae^+r^+\underline Q_1\Lambda
_+^a\Lambda _+^{1-a}e^+ K_{\underline
Q_1^+}\underline{\varrho _i\varphi }+
r^+\underline \zeta _i^2\Lambda _-^ae^+r^+\underline Q_1\Lambda
_+^a(1-\underline\zeta _i^1)\Lambda _+^{1-a}e^+ K_{\underline
Q_1^+}\underline{\varrho _i\varphi }.
\endaligned
$$
The first term equals
$$
r^+\underline \zeta _i^2\Lambda _-^ae^+r^+\underline Q_1\Lambda _+^{1}e^+ K_{\underline
Q_1^+}\underline{\varrho _i\varphi },\tag5.32
$$
and will be treated further below. The second term equals
$$
\aligned
r^+\underline \zeta _i^2&\Lambda _-^ae^+r^+\underline Q_1\Lambda
_+^a(1-\underline\zeta _i^1)\Lambda _+^{1-a}e^+ r^+\Lambda
_+^{-1}\widetilde{\underline Q}_1^+\underline\zeta _i^0(\underline{\varrho _i\varphi }\otimes
\delta )\\
&=r^+\underline \zeta _i^2\Lambda _-^ae^+r^+\underline Q_1\Lambda
_+^a(1-\underline\zeta _i^1)\Lambda _+^{1-a}\Lambda
_+^{-1}\widetilde{\underline Q}_1^+\underline\zeta _i^0(\underline{\varrho _i\varphi }\otimes
\delta )\\
&=r^+\underline \zeta _i^2\Lambda _-^ae^+r^+\underline Q_1\Lambda
_+^a(1-\underline\zeta _i^1)\Lambda _+^{-a}\widetilde{\underline Q}_1^+\underline\zeta _i^0(\underline{\varrho _i\varphi }\otimes
\delta ),
\endaligned
$$
where $\Lambda
_+^a(1-\underline\zeta _i^1)\Lambda _+^{-a}$ is a $\psi $do of order 0
having the 0-transmission property; we used in the proof that $(I-e^+r^+)\Lambda
_+^{-1}\widetilde{\underline Q}_1^+\underline\zeta _i^0(\underline{\varrho _i\varphi }\otimes
\delta )=0$. In this term,
% is further transformed as
$$
\underline{\varrho _i\varphi }\mapsto
r^+\underline Q_1\Lambda
_+^a(1-\underline\zeta _i^1)\Lambda _+^{-a}\widetilde{\underline Q}_1^+\underline\zeta _i^0(\underline{\varrho _i\varphi }\otimes \delta )\tag5.33
$$
is a Poisson operator of order 0, derived from a generalized $\psi $do
whose symbol contains the factors $1-\underline\zeta _i^1$ and
$\underline\zeta _i^0$ with disjoint supports; hence it must be a Poisson operator of order $-\infty
$. Then it maps into $\bigcap _t\ol H^t(\rnp)$, and the composition with
$r^+\underline \zeta _i^2\Lambda _-^ae^+$ maps into $r^+C_0^\infty
(V_i)$. In conclusion, this term reduces to a smoothing contribution $\Cal
S_{i,3}\varphi $.

We now continue with the term in (5.32). Here 
$$
\aligned
r^+\underline Q_1\Lambda _+^{1}e^+ K_{\underline
Q_1^+}\underline{\varrho _i\varphi }&=r^+\underline Q_1\Lambda _+^{1}e^+r^+
\Lambda _+^{-1}\widetilde {\underline Q}_1^+ (\underline{\varrho
_i\varphi }\otimes\delta )\\
&=r^+(\underline Q_1^-\underline Q_1^++\Cal R_1)
\Lambda _+^{1}
\Lambda _+^{-1}\widetilde {\underline Q}_1^+ (\underline{\varrho
_i\varphi }\otimes\delta )\\
&=r^+(\underline Q_1^-\underline Q_1^++\Cal R_1)
\widetilde {\underline Q}_1^+ (\underline{\varrho
_i\varphi }\otimes\delta ),
\endaligned
$$
where $\Cal R_1$ has symbol in $S^{-\infty }(\Cal H_{-1})$.
The contribution from $\Cal R_1$ is seen as above to be smoothing, and
$$
r^+\underline Q_1^-\underline Q_1^+
\widetilde {\underline Q}_1^+ (\underline{\varrho
_i\varphi }\otimes\delta )=
r^+\underline Q_1^-(I+\Cal R_2) (\underline{\varrho
_i\varphi }\otimes\delta ),
$$
where the contribution from $\Cal R_2$ is likewise seen to be
smoothing. Finally, 
$$
r^+\underline Q_1^- (\underline{\varrho
_i\varphi }\otimes\delta )=0,
$$
since $\underline Q_1^-$ preserves support in $\crnm$. Collecting the
contributions, carried back to the coordinate patches $U_i$, we find that
$$
r^+PK_D\varphi =\Cal S\varphi ,\tag5.34
$$
where $\Cal S$ maps into $C^\infty (\comega)$.

We also have to show that $\gamma _0^{a-1}K_D=I$. For functions $u$ on
$\Omega $, $\gamma _0^{a-1}u=\Gamma (a)\gamma _0(d^{1-a}u)$, where
$d(x)=\operatorname{dist}(x,\partial\Omega )$. With our special choice
of local coordinates, this carries over from $U_i\cap\Omega $ to
$V_i\cap\rnp$ as $\Gamma (a)\gamma _0(x_n^{1-a}\underline u)$ (the trace at $x_n=0$), when
$\operatorname{supp}u\subset U_i$.
 When $u$ is
multiplied by a function $\zeta $, $\gamma _0^{a-1}(\zeta u)=\gamma
_0(\zeta )\gamma _0^{a-1}(u)$, which carries over as $\gamma
_0(\underline\zeta )\gamma _0^{a-1}(\underline u)$. Recall from
(5.25) (or (A.16)) that the boundary value
$\gamma _0^{a-1}(\underline u)$ can
also be described as $\gamma _0(\Xi _+^{a-1}\underline u)$ or $\gamma _0(\Lambda
_+^{a-1}\underline u)$. 

Let $\varphi $ be given as above. Write 
$\varphi
=\sum_{0\le i\le I_0}\varrho _i\varphi $, then
$$
\gamma _0^{a-1}K_D\varphi =\sum_i\gamma _0^{a-1}(\underline\zeta _i^1\underline K_D^{(i)}\underline{\varrho _i\varphi } )\circ \kappa _i.
$$
Let us consider the $i$'th
 piece, carried over to $V_i$:
%, allowing $(i)$ to be omitted from the notation:
$$
\aligned
\gamma _0^{a-1}(\underline\zeta _i^1\underline
K_D^{(i)}\underline{\varrho _i\varphi } )&=
\gamma _0(\underline\zeta _i^{1})\gamma _0^{a-1}(\underline
K_D^{(i)}\underline{\varrho _i\varphi } )\\
&
=\gamma _0(\underline\zeta _i^{1})
\gamma _0^{a-1}(\Lambda _+^{1-a}e^+ r^+\Lambda
_+^{-1}\widetilde{\underline Q}_1^{(i)+}(\underline{\varrho _i\varphi }\otimes
\delta ))\\
&
=\gamma _0(\underline\zeta _i^{1})
\gamma _{0}(r^+\Lambda _+^{-1}\widetilde{\underline Q}_1^{(i)+}(\underline{\varrho _i\varphi }\otimes
\delta ))\\
&
=\gamma _0(\underline\zeta _i^{1})
        \underline{\varrho _i\varphi }=\underline{\varrho _i\varphi }.
\endaligned
$$
In the step leading to the last line, we used that $\widetilde{\underline Q}_1^{(i)+}$ has a symbol $1+f^+$, where $f^+\in
S^0(\Cal H^+)$ and hence does not contribute to the boundary value, so
only $K_0$ remains (as in
the proof of Theorem 5.2).

Carrying the formulas back to the $U_i$ and
summing over $i$, we obtain the conclusion 
$$
\gamma _0^{a-1}K_D\varphi =\varphi .\tag5.35
$$

The last part of the proof goes in the same way as in Corollary 5.3.\qed

\enddemo

\example{Remark 5.5}
In the above proof, one can replace $K_D$ by the operator defined
using the family $\Xi _\pm^t$ instead of $\Lambda _\pm^t$ in the
local coordinate systems; it then gets the form
$$
K_D\varphi =\sum_{i=1}^{I_0}(\underline\zeta _i^1\Xi _+^{1-a}e^+
K^{(i)}_{{\underline Q^+} }\underline{\varrho _i\varphi })\circ \kappa
_i,\text{ where } K^{(i)}_{{\underline Q^+} }\underline{\varrho _i\varphi }=r^+\Xi _+^{-1}\widetilde {\underline Q} ^{(i)+}(\underline{\varrho _i\varphi }
\otimes\delta (x_n)).\tag5.36
$$ 
The calculations go as above, except that in the treatment of the term
corresponding to (5.33),
$\Xi 
_+^a(1-\underline\zeta _i^1)\Xi  _+^{-a} $ is only a generalized $\psi
$do; but it can be checked to have a symbol of the form of a function
plus a term in $S^0(\Cal H^+)$ (only integer powers of $\Xi _+^1$
appear in the terms calculated by the Leibniz product formula), and the Poisson operator construction
goes through.
\endexample

There is a certain ``uniqueness modulo smoothing operators'':

\proclaim{Corollary 5.6} In the situation of Theorem {\rm 5.4}, if
$\,^t\pmatrix r^+P &\gamma
_0^{a-1}\endpmatrix$ has the inverse \linebreak $\pmatrix R_D&K_D\endpmatrix$, 
 and $K_{D,1}$ is a Poisson-like operator
constructed as in Theorem {\rm 5.4} or Remark {\rm 5.5}, then 
$$
K_D-K_{D,1}\colon H^{s-a+\frac12}(\partial\Omega )\to \Cal E_{a}(\comega).\tag5.37
$$

\endproclaim

\demo{Proof} We have that
$$
\aligned
\pmatrix 0&K_D-K_{D,1}\endpmatrix &=\pmatrix R_D&K_D\endpmatrix - \pmatrix R_D&K_{D,1}\endpmatrix\\
&= \pmatrix R_D&K_D\endpmatrix \bigl[\pmatrix r^+P\\\gamma
_0^{a-1}\endpmatrix \pmatrix R_D&K_{D,1}\endpmatrix- \pmatrix \Cal S_1&\Cal
S\\0&0\endpmatrix\bigr]\\
&\quad- \pmatrix R_D&K_D\endpmatrix \pmatrix
r^+P\\\gamma _0^{a-1}\endpmatrix \pmatrix R_D&K_{D,1}\endpmatrix\\
&=\pmatrix R_D\Cal S_1& R_D\Cal S\endpmatrix.
\endaligned
$$
Since $R_D\colon C^\infty (\comega)\to \Cal E_{a} (\comega )$,  $K_D-K_{D,1}=R_D\Cal
S$ satisfies (5.37).\qed

\enddemo

\subhead 6. The
Dirichlet-to-Neumann operator, elliptic Neumann problems
\endsubhead 

We can now calculate the {\it Dirichlet-to-Neumann} operator $$
S_{DN}=\gamma _1^{a-1}K_D,\tag6.1$$
finding its symbol in local coordinates.
In particular, we determine its principal symbol, which gives a
criterion for ellipticity of the Neumann problem.

\proclaim{Theorem 6.1} 
$1^\circ$ For $P$ considered on $\crnp$ as in Theorem {\rm 5.2}, with $K_D$ defined by {\rm
(5.16--19)},  $S_{DN}$ is the first-order $\psi $do
$S_{DN}=\OP'(s_{DN}(x',\xi '))$ in ${\Bbb R}^{n-1}$
with symbol (cf.\ {\rm (5.10))}
$$
\aligned 
s_{DN}(x',\xi ')&=-s_{ Q ^+}(x',\xi ') -a\ang{\xi
'},\text{ where }\\
s_{ Q ^{+}}(x',\xi ')&=\tfrac1{2\pi }\int^+ h^+ q^{+}(x',0,\xi
)\,d\xi _n =\lim_{z_n\to 0+}\Cal
F^{-1}_{\xi _n\to z_n} h^+q^+(x',0,\xi ).
\endaligned\tag6.2
$$
In particular, if $P=(1-\Delta )^a$, then $S_{DN}=-a\ang{D'}$.

$2^\circ$ For $P$ considered on $\comega$ as in Theorem {\rm 5.4}, with $K_D$ defined by {\rm
(5.36)},  $S_{DN}$ is the first-order $\psi $do on $\partial\Omega $
described by
$$
S_{DN}\varphi =\sum_{i=1}^{I_0}(\gamma _0(\underline\zeta
_i^1)\operatorname{OP}'(-s_{\underline{
Q} ^{(i)+}}(x',\xi')-a\ang{\xi '})\underline{\varrho _i\varphi })\circ \kappa _i.\tag6.3
$$
where $s_{\underline{
Q} ^{(i)+}}(x',\xi')$ is constructed from the symbol
of $\underline{
Q} ^{(i)+}$ as in {\rm (6.2)}.
\endproclaim 

\demo{Proof} $1^\circ$. Recall  from (3.8) that $\gamma _1^{a-1}u=\gamma
_0\partial_n\Xi _+^{a-1}u-(a-1)\ang{D'}\gamma _0^{a-1}u$. We have
immediately for the second term:
$$
-(a-1)\ang{D'}\gamma _0^{a-1}K_D=-(a-1)\ang{D'},\tag6.4
$$
by (5.21). For the first term we apply composition rules from the
Boutet de Monvel calculus:
$$
\aligned
\gamma _0\partial_n\Xi _+^{a-1}K_D&=\gamma
_0\partial_n\Xi _+^{a-1}\Xi _+^{1-a}e^+\operatorname{OPK}(\chi
_+^{-1}(\xi )\#\widetilde q^{+\prime} (x',0,\xi ))\\
&=\gamma
_0\operatorname{OPK}(i\xi _n(\ang{\xi '}+i\xi _{n})^{-1}\#\widetilde q^{+\prime} (x',0,\xi )).\endaligned
$$
This is the $\psi $do with symbol
$$
%\aligned
\tfrac1{2\pi }\int^+\frac{i\xi _n}{\ang{\xi '}+i\xi _n}\#\tilde
q^{+\prime} \,d\xi _n %&
=\tfrac1{2\pi }\int^+(1-\frac{\ang{\xi
'}}{\ang{\xi '}+i\xi _n})\#(1+h^+\widetilde q^{+\prime} )\,d\xi _n%\\
%&
=-\ang{\xi '}+\tfrac1{2\pi }\int^+h^+\widetilde q^{+\prime} \,d\xi _n.%\\
%&=-\ang{\xi '}+s_{\widetilde Q ^{+\prime}}(x',\xi '),
%\endaligned
$$
%cf.\ (6.2). 
We have here used that $\int^+1\,d\xi _n=0$, $\frac1{2\pi }\int^+({\ang{\xi
'}+i\xi _n})^{-1}\,d\xi _n=1$, and that the plus-integral
vanishes on functions in $\Cal H^+$ that are $O(\xi _n^{-2})$,
hence on $(\ang{\xi '}+i\xi _n)^{-1}\#h^+\widetilde q^{+\prime} $.

The two terms together give the $\psi $do with symbol
$$
-\ang{\xi '}+\tfrac1{2\pi }\int^+h^+\widetilde q^{+\prime} \,d\xi _n-(a-1)\ang{\xi '}=\tfrac1{2\pi }\int^+h^+\widetilde q^{+\prime} \,d\xi _n-a\ang{\xi '}.
$$
%showing (6.2).

There is a further reduction of the plus-integral. We know that
$q^+(x,\xi )$ is of the form $q^+=1+f^+$, where $f^+=h^+q^+$ lies in
$S^0(\Cal H^+)$. The parametrix $\widetilde q^+$ of $q^+$ has the expansion
$$
\widetilde q^+\sim 1-f^++f^+\# f^+ -\dots +(-1)^k(f^+)^{\#k}+\dots,
$$
where $(-1)^k(f^+)^{\#k}\in S^0(\Cal H^+\cap \Cal H_{-2})$ for $k\ge
2$ (all terms have at least two factors in $\Cal H^+$), so we can
assume
that
$\widetilde q^+$ is of the form
$$
\widetilde q^+=1-f^++r,\quad r\in S^0(\Cal H^+\cap \Cal H_{-2}).
$$
Then 
$$
\int^+h^+\widetilde q^+\,d\xi _n=\int^+\widetilde q^+\,d\xi _n=
\int^+(1-f^++r)\,d\xi _n=-\int^+f^+\,d\xi _n=-\int^+h^+q^+\,d\xi _n.
$$
Now it is actually the $(x',y_n)$-form  of the symbol of $\widetilde
Q^+$ that is used instead of the $x$-form, cf.\ (A.11), but this does not
change the value, as we shall now show: 

The preceding calculations are
true also for $\widetilde q^{+\prime}$  and $q^{+\prime}$, so we
arrive at having to calculate  $\frac1{2\pi }\int^+h^+q^{+\prime}(x',0,\xi )\,d\xi _n$.
As in the proof details after (4.15), we can use the observation from \cite{G90} Lemma 10.18 that for a function $\varphi (\xi
_n)\in \Cal H^+$,
$
\tfrac1{2\pi }\int^+\varphi (x_n)\,d\xi _n=\lim_{z_n\to 0+}[\Cal
F^{-1}_{\xi _n\to z_n}\varphi ](z_n)$.
Since
$$
[\Cal
F^{-1}_{\xi _n\to z_n} {h^+q^{+\prime}}(x',0,\xi )](z_n)\sim \sum_{j\in {\Bbb N}_0}\tfrac
1{j!}z_n^j[\Cal
F^{-1}_{\xi _n\to z_n} h^+\partial_{x_n}^jq^+(x',0,\xi )](z_n),
$$
the limits for $z_n\to 0+$ satisfy 
$$
\lim_{z_n\to 0+}[\Cal
F^{-1}_{\xi _n\to z_n} {h^+q^{+\prime}}(x',0,\xi )](z_n)=\lim_{z_n\to 0+}[\Cal
F^{-1}_{\xi _n\to z_n} h^+q^+(x',0,\xi )](z_n).
$$
Thus $h^+q^{+\prime}$
gives the same value as $h^+q^+$ in the plus-integral, and formula
(6.2) follows.

When $P=(1-\Delta )^a$, then $Q=I$ with $h^+\widetilde q^+=0$, so only the term $-a\ang{ \xi '}$
remains. The formula was shown in this case in \cite{G14}, Appendix.

$2^\circ$. We can assume that the $\underline \zeta _i^1$ are constant
in $x_n$ for small $x_n$, then
$$
\gamma _1^{a-1}\underline u=\Gamma (a+1)\gamma _0(\partial_n(\underline\zeta
_i^1\underline u/x_n^{a-1}))
=\Gamma (a+1)\gamma _0(\underline\zeta _i^1)\gamma _0(\partial_n(\underline u/x_n^{a-1}))=
\gamma _0(\underline\zeta _i^1)\gamma _1^{a-1}(\underline u).
$$
Now 
$$
\gamma _1^{a-1}K_D\varphi =\sum_i(\gamma _0(\underline\zeta _i^1)\gamma _1^{a-1}(\Xi _+^{1-a}e^+ K^{(i)}_{{\underline Q^+} }\underline{\varrho _i\varphi }))\circ \kappa _i,
\tag6.5
$$
where each term is calculated as under $1^\circ$. This leads to (6.3).\qed
\enddemo

\proclaim{Theorem 6.2} Hypotheses as in Theorem {\rm 6.1}. 

The principal symbol of $S_{DN}$ is 
$$
s_{DN,0}(x',\xi ')=-\tfrac1{2\pi }\int^+h^+q^+_{0}(x',0,\xi )\,d\xi _n-a|\xi '|,\tag6.6
$$
where $q^+_{0}$ is the 
 plus-factor in the principal symbol $q_0(x,\xi )=s_0^{-1}p_0(x,\xi )|\xi |^{-2a}$, constructed as in \cite{G16},
Th.\ {\rm 2.6}, in the local coordinates used in Theorem {\rm 6.1}.
Here
$$
\tfrac1{2\pi }\int^+h^+q^+_{0}(x',0,\xi )\,d\xi _n=\lim_{z_n\to 0+}\Cal
F^{-1}_{\xi _n\to z_n}\log q_0(x',0,\xi ).\tag6.7
$$
%$q_0(x,\xi )=p_0(x,\xi )|\xi | ^{-2a}$.
%the boundary value for $x_n\searrow 0$.

The Neumann problem defined for $u\in H^{(a-1)(s)}(\comega)$, $s>a+\frac12$,
$$
r^+Pu=f\text{ on }\Omega ,\quad  \gamma _1^{a-1}u=\psi ,\tag 6.8
$$
has a parametrix $$
\pmatrix R_N& K_N\endpmatrix\colon \ol H^{s-2a}(\Omega )\times
H^{s-a-\frac12}(\partial\Omega )\to  H^{(a-1)(s)}(\comega),\tag6.9
$$
if and only $s_{DN}(x',\xi ')$ is nonvanishing for $\xi '\ne 0$, i.e., $S_{DN}$ is elliptic. In the affirmative case, a
parametrix is
$$
\pmatrix R_N& K_N\endpmatrix=\pmatrix (I-K_D\widetilde S_{DN}\gamma _1^{a-1})R_D& K_D\widetilde S_{DN}\endpmatrix,\tag6.10
$$
where $\widetilde S_{DN}$ is a parametrix of $S_{DN}$. The Neumann problem
 is then said to be elliptic.
\endproclaim

\demo{Proof} It is known that the principal symbol of a $\psi $do on a
manifold (in this case $\partial\Omega $) is an invariant function of
the cotangent variables $(x',\xi ')$, so we just need to indicate a
way to find it; this goes via the localization (6.3), when we moreover
use that $\zeta _i^1\varrho _i=\varrho  _i$ for each $i$. Then the formula
$$
s_{DN,0}(x',\xi ')=-\tfrac1{2\pi }\int^+h^+ q^+_{0}(x',0,\xi )\,d\xi _n-a|\xi '|,
$$
follows from Theorem 6.1.

For the description in (6.7), we recall from the proof of Theorem 2.6 in \cite{G16} that $q^+_0=\exp
\psi _+$, where $\psi _+=h^+\psi $ with  $\psi =\log q_0$.
Here $q_0^+$ 
has the expansion
$$
q_0^+=1+\psi _++\sum_{k\ge 2}
\tfrac1{k!} \psi _+^k,
$$
where the sum over $k\ge 2$ is $O(\ang{\xi _n}^{-2})$, so this, as well as the term 1, gives 0
when plus-integrals are calculated, and only $\psi _+$ %resp.\ $-\psi_+$ 
remains.
 Thus $$
\tfrac1{2\pi
}\int^+h^+ q^+_0\,dx_n=
\tfrac1{2\pi
}\int^+\psi _+\,dx_n=
\lim_{z_n\to 0+}\Cal F^{-1}_{\xi _n\to z_n} \psi_+=\lim_{z_n\to 0+}\Cal F^{-1}_{\xi _n\to z_n} \psi,$$
since $\psi _-=\psi -\psi _+=h^-\psi $ has $\Cal F^{-1}_{\xi _n\to z_n} \psi
_-$ supported in $\crm$.

For the second statement, it is easily checked that if $S_{DN}$ is
elliptic with a parametrix $\widetilde S_{DN}$, then a parametrix
for (6.8) can be constructed from the parametrix for the
Dirichlet problem by the formula (6.9). (Similar considerations are
carried out in \cite{G14}, Sect.\ 4B.) Conversely, if the Neumann problem
has a parametrix with the asserted mapping property, then there is a
Fredholm mapping
from $\varphi _1 =\gamma _1^{a-1}u $ to $\varphi _0= \gamma _0^{a-1}u $ for
solutions $u$
to the Neumann resp.\ Dirichlet problem with the same value of
$r^+Pu$. This gives a parametrix of 
$S_{DN}$; hence it is elliptic.  
\qed
\enddemo

\example{Example 6.3} To illustrate some of the formulas, let
us consider a very simple example with an easy factorization. Let $P=(\beta
|D'|^2+D_n^2+\lambda )^a$ with $\beta >0$, $\lambda \ge0$, in the situation $\rnp\subset {\Bbb R}^n$. Then 
$$
Q =(\beta
|D'|^2+D_n^2+\lambda )^a(1-\Delta )^{-a}=\OP(q^-q^+), \quad q^\pm =
\bigl(\frac{(\beta |\xi '|^2+\lambda )^\frac12\pm i\xi _n}{\ang{\xi '}\pm i\xi _n}\bigr)^{a}.
$$
Here $ q^+$
%=(q^+)^{-1}$ 
has the expansion
$$
\aligned
q^+&=\bigl(\frac{(\beta |\xi '|^2+\lambda )^{\frac12}+ i\xi _n}{\ang{\xi '}+ i\xi
_n}\bigr)^{a}= \bigl(1+\frac{(\beta |\xi '|^2+\lambda )^{\frac12} -\ang{\xi '}}{\ang{\xi '}+ i\xi
_n}\bigr)^{a}\\
&=1+a\frac{(\beta |\xi '|^2+\lambda )^{\frac12} -\ang{\xi '}}{\ang{\xi '}+ i\xi
_n}+O(\ang{\xi _n}^{-2}),
\endaligned
$$
 so 
$$
\aligned
-\tfrac1{2\pi }\int^+h^+ q^+\,d\xi _n&=-a((\beta |\xi '|^2+\lambda )^\frac12-\ang{\xi '})\tfrac1{2\pi }\int^+\frac1{\ang{\xi '}+ i\xi
_n}\,d\xi _n\\
&=-a(\beta |\xi '|^2+\lambda )^\frac12+a\ang{\xi '}.\endaligned
$$
Then 
$$
s_{DN}=-a(\beta |\xi '|^2+\lambda )^\frac12+a\ang{\xi '}-a\ang{\xi '}=
-a(\beta |\xi '|^2+\lambda )^\frac12,\quad
s_{DN,0}=-a\beta ^\frac12|\xi '|.\tag6.10
$$
In this case $s_{DN}$ is elliptic for all choices of $\beta >0, 
\lambda \ge0$, and the
Neumann problems for these operators are elliptic. 
If $\lambda $ is replaced by a nonnegative function $V(x)$, the
calculations are valid on the principal symbol level.

\endexample

In the paper \cite{G15} the continuity properties of $\,^t\pmatrix r^+P&
\gamma _0^{a-1}\endpmatrix$ and its parametrix $\pmatrix R_D&
K_D\endpmatrix$ were shown in $H^s_p$-Sobolev spaces (essentially, see
Ths.\ 4.2, 5.1 and 6.5 there), and in \cite{G14}, which was written
after \cite{G15}, they were extended to large families of Besov and
Triebel-Lizorkin spaces, including results for   elliptic Neumann
problems (see Ths.\ 3.2, 3.5 and 4.3 there). Only the structure of $K_D$ and
$K_N$ in the case of a curved domain was not explained in detail.
 We now have the full explanation above of how $K_D$ and $K_N$ 
consist of operators from the Boutet
de Monvel calculus composed with operators $\Xi _+^{1-a}$ or $\Lambda _+
^{1-a}$ in local coordinates, so the mapping properties extend readily
from the $H^s_2$-scales
to the $H^s_p$-scales as in \cite{G15}, and to the Besov and Triebel-Lizorkin scales as
in \cite{G14} Sect.\ 3 (based on Johnsen
\cite{J96} and its references). 
We can therefore conclude from \cite{G15, G14} the following formulations of
the mapping properties of the parametrices, for smooth bounded domains
$\Omega $. (We recall that $H^s_p=F^s_{p,2}$ and $C^s_*=B^s_{\infty ,\infty }$.)

\proclaim{Corollary 6.4} Let $1\le p\le \infty $ and $0<q\le \infty $
(with $p<\infty $ in the $F$-cases,  {\rm (6.11), (6.13)} and the first
line of {\rm (6.15)} below). 

The parametrix of the Dirichlet problem maps
as follows, for $s>a-1/p'$:
$$\aligned
\pmatrix R_D& K_D\endpmatrix&\colon\ol
F_{p,q}^{s-2a}(\Omega )\times
B_{p,p}^{s-a+1/p'}(\partial\Omega )\to
F_{p,q}^{(a-1)(s)}(\comega),\text{ in particular}\\ 
\pmatrix R_D& K_D\endpmatrix
&\colon \ol H_p^{s-2a}(\Omega )\times
B_{p,p}^{s-a+1/{p'}}(\partial\Omega )\to H_p^{(a-1)(s)}(\comega).
\endaligned\tag6.11$$
Moreover, for $s>a-1/p'$, $t>a-1$,
$$\aligned
\pmatrix R_D& K_D\endpmatrix&\colon\ol
B_{p,q}^{s-2a}(\Omega )\times
B_{p,q}^{s-a+1/p'}(\partial\Omega )\to B_{p,q}^{(a-1)(s)}(\comega),
\text{ in particular}\\
\pmatrix R_D& K_D\endpmatrix &\colon 
\ol
C_*^{t-2a}(\Omega )\times
C_*^{t-a+1}(\partial\Omega )\to C_*^{(a-1)(t)}(\comega).
\endaligned\tag6.12$$

The parametrix of the Neumann problem (when elliptic) maps
as follows, for $s>a+1/p$:
$$\aligned
\pmatrix R_N& K_N\endpmatrix&\colon\ol
F_{p,q}^{s-2a}(\Omega )\times
B_{p,p}^{s-a-1/p}(\partial\Omega )\to F_{p,q}^{(a-1)(s)}(\comega),\\ 
\pmatrix R_N& K_N\endpmatrix
&\colon \ol H_p^{s-2a}(\Omega )\times
B_{p,p}^{s-a-1/{p}}(\partial\Omega )\to H_p^{(a-1)(s)}(\comega).
\endaligned\tag6.13$$
Moreover, for $s>a+1/p$, $t>a$,
$$\aligned
\pmatrix R_N& K_N\endpmatrix&\colon\ol
B_{p,q}^{s-2a}(\Omega )\times
B_{p,q}^{s-a-1/p}(\partial\Omega )\to B_{p,q}^{(a-1)(s)}(\comega),\\
\pmatrix R_N& K_N\endpmatrix &\colon 
\ol
C_*^{t-2a}(\Omega )\times
C_*^{t-a}(\partial\Omega )\to C_*^{(a-1)(t)}(\comega). 
\endaligned\tag6.14$$

Being a $\psi $do of order $1$, the Dirichlet-to-Neumann operator
$S_{DN}$ maps as follows, for all $s\in{\Bbb R}$:
$$\aligned
S_{DN}&\colon F_{p,q}^{s}(\partial\Omega )\to
F_{p,q}^{s-1}(\partial\Omega ),\quad
S_{DN}\colon H_{p}^{s}(\partial\Omega )\to H_{p}^{s-1}(\partial\Omega
),\\
S_{DN}&\colon B_{p,q}^{s}(\partial\Omega )\to
B_{p,q}^{s-1}(\partial\Omega ),\quad
S_{DN}\colon C_*^{s}(\partial\Omega )\to C_*^{s-1}(\partial\Omega ).
\endaligned\tag6.15
$$
\endproclaim

\example{Remark 6.5} In connection with the Neumann trace operator $\gamma
_1^{a-1}$ there are many other meaningful boundary conditions for
$r^+P$, namely conditions  of the (Robin-like) form
$$
\gamma _1^{a-1}u+L\gamma _0^{a-1}u=\psi ,\tag6.16
$$
where $L$ is a first-order $\psi $do on $\partial\Omega $; they are
{\it local} when $L$ is a differential operator. The problem for $r^+P$ with
condition (6.16) is {\it elliptic} when $S_{DN}+L$ is elliptic, i.e.,
in local coordinates, $s_{DN,0}(x',\xi ')+l_0(x',\xi ')\ne 0$ for $\xi
'\ne 0$. Then the problem is Fredholm solvable in scales of spaces as
in Theorem 6.2 and Corollary 6.4.

In particular, when $P$ is principally like $(-\Delta )^a$ (i.e.,
$P-(-\Delta )^a$ is of order $2a-1$),  we have, when $L$ in
local  coordinates at the boundary has the principal part $\sum
_{j=1}^{n-1}b_j(x')\partial_j $, that
ellipticity holds when $-a|\xi '|+b(x')\cdot i \xi '$ is
nonvanishing for $\xi '\ne 0$. This is satisfied when the vector $b$ is real; for complex $b$ it holds
when $|\operatorname{Im}b(x')|<a$ for all $x'$.

\endexample

Finally, some remarks on nonsmooth situations: Abels \cite{A05a, A05b}
has established a nonsmooth version of results from the Boutet de
Monvel calculus in $H^s_p$ and $B^s_{p,q}$-spaces for operators with 
for example H\"older continuous $x$-dependence. This is based on earlier works
on nonsmooth $\psi $do's, in particular Marschall
\cite{M88}. In this generalization, the range of possible $s$ is
limited  to
specific bounded intervals. Certainly, much of the above can be
extended to such operators, when $\Omega $ is smooth. To allow
nonsmoothness of $\Omega $, one would need to handle nonsmooth
coordinate changes, but here the results are scarce. --- Another strategy
would be to use the results known in smooth cases in  combination with
perturbation arguments (for example, Vishik and Eskin \cite{E81} used such techniques to pass
from constant-coefficient cases to variable coefficients). Much
remains to be investigated for such questions.

The results from the smooth studies can serve as a model for what one
would want to show in nonsmooth cases. In comparison, many results
through the years involving the Laplacian $-\Delta $ are
generalizations  of old and well-established results in smooth
cases, for example  based on Green's formula. The ``old, well-established'' results have not been available
in the case of the fractional Laplacian  $(-\Delta )^a$; the present
strivings with $\psi $do methods help to fill that gap.

\subhead Appendix
\endsubhead

We here collect some notation and results from the Boutet de Monvel
calculus (as exposed in \cite{B71, G90, G96, S01, G09}) and from  \cite{G16}, that are used
in the text.

The operators below are defined for the situation $\rnp\subset{\Bbb
R}^n$. Besides pseudodifferential operators ($\psi $do's) on
${\Bbb R}^n$, cf.\ (2.2), there are {\it Poisson operators} $K$ from
${\Bbb R}^{n-1}$ to $\rnp$, and {\it trace
operators} $T$ of class 0 from $\rnp$ to ${\Bbb R}^{n-1}$:
$$\aligned
(K\varphi )(x)&=\operatorname{OPK}(k(x',\xi ))\varphi =(2\pi )^{-n}\int e^{ix\cdot \xi }k(x',\xi )\hat\varphi (\xi ')\, d\xi \\
(Tu)(x')&=\operatorname{OPT}(t(x',\xi ))u=(2\pi )^{-n}\int e^{ix'\cdot \xi '}t(x',\xi )\widehat {e^+u}\, d\xi ,
\endaligned\tag A.1
$$
with suitable interpretations of the integrals, and $\psi $do's $S=\OP'(s(x',\xi '))$ on ${\Bbb R}^{n-1}$ defined as
in (2.2). In general, the trace operators of class 0 are supplied with standard
trace operators $\sum_{0\le j<r}S_j\gamma _j$ with $\psi $do
coefficients $S_j$; then they are of class $r$. (We shall
not need {\it singular Green operators} in the present paper.) The
symbols in (A.1) behave w.r.t.\ $\xi _n$ as elements of the spaces
$\Cal H^+$ resp.\ $\Cal H^-_{-1}$ of functions $f(t)$ that we shall
now recall:

For $d\in{\Bbb Z}$, 
$\Cal H_{d}$ denotes the space of $C^\infty $-functions $f(t)$
on ${\Bbb R}$ such that $k(\tau )=\tau ^df(\tau ^{-1})$ coincides with
a $C^\infty $-function for $-1<\tau <1$ (i.e., the
derivatives of $f$ match in a good way for $t\to\pm\infty $). These
spaces are Fourier transforms of spaces of distributions on ${\Bbb R}$:
$$
\Cal H_{-1}=\Cal F(e^-\Cal S_-\oplus e^+\Cal S_+),\quad \Cal
H_{d}=\Cal H_{-1}\oplus {\Bbb C}_d[t] \text{ for }d\ge 0,\tag A.2
$$
where $\Cal S_{\pm}=r^\pm \Cal S({\Bbb R})=\Cal S(\crpm)$, and
${\Bbb C}_d[t]$ stands for the space of polynomials of degree $\le d$
in $t$. We also denote
$$
\Cal H=\bigcup_d \Cal H_d=\Cal H_{-1}\oplus {\Bbb C}[t]
, \quad {\Bbb C}[t]=\bigcup_d {\Bbb C}_d[t].$$
With a slight asymmetry, one defines 
$$
\aligned
\Cal H^+&=\Cal F (e^+\Cal S_+), \; \Cal H^-=\Cal F (e^-\Cal
S_-)\oplus {\Bbb C}[t],\; \Cal H^-_{d}=\Cal F (e^-\Cal
S_-)\oplus {\Bbb C}_d[t] \text{ for }d\ge -1,\text{ so}
\\
\Cal H&=\Cal H^+\oplus\Cal H^-=\Cal H^+\oplus\Cal H^-_{-1}\oplus {\Bbb C}[t],
\text{ with }
\Cal H_{-1}=\Cal
H^+\oplus\Cal H^-_{-1};
\endaligned
  \tag A.3
$$
more generally $ \Cal H_d=\Cal H^+\oplus\Cal H^-_d$ for $d\ge -1$.
The mappings $h^\pm$ are defined on $\Cal H$
such that they are complementing projections  with
ranges $\Cal H^\pm$:
$$%\aligned
h^+\Cal H=\Cal H^+,\quad h^-\Cal H=\Cal H^- ,\text{ in particular, } h^-\Cal H_d=\Cal H^-_d 
\text{ for
}d\ge -1.\tag A.4
%\endaligned
$$
Note that $h^+$ and $h^-$ are essentially the Fourier transforms of
the projections $e^+r^+$ and $e^-r^-$ for functions on ${\Bbb R}$;
$\F^{-1}h^-\F$ moreover preserves distributions supported in $\{0\}$.

Furthermore, $h_{-1}$ denotes the projection from $\Cal H$ to $\Cal
H_{-1}$ that removes the polynomial part.
The space $\Cal H^-_{-1}$ equals the space of conjugates of functions
in $\Cal H^+$. $\Cal H^+$ can also be denoted
$\Cal H^+_{-1}$ when convenient. When $f\in\Cal H_{-1}$,
$\overline{h^-f}=h^+(\overline f)$.

The symbol spaces for Poisson resp.\ trace operators of class 0 are the spaces
$S^m(\Cal H^+)$ resp.\ $S^m(\Cal H^-_{-1})$. Here $S^m(\Cal
H^+)=S_{1,0}^m({\Bbb R}^{m},{\Bbb R}^{n-1},\Cal H^+)$ consists of
functions $f(X,\xi ',\xi _n)$ that are in $\Cal H^+$ w.r.t.\ $\xi _n$
and satisfy the estimates
$$
\| D_{X}^\beta D_{\xi '}^\alpha D_{\xi _n}^k h_{-1}(\xi _n^{k'}f(X,\xi
',\xi _n))\|_{L_2({\Bbb R})}\le C_{\alpha ,\beta ,k,k'} \ang{\xi
'}^{m+\frac12-k+k'-|\alpha |},\tag A.5 
$$
for all indices. The symbols considered here moreover have asymptotic expansions in terms
$f_j$ that are homogeneous in $(\xi ',\xi _n)$ of degree $m-j$ for $|\xi
'|\ge 1$, such that $f-\sum_{0\le j<J}f_j$ is in $S^{m-J}(\Cal H^+)$
for all $J$. Definitions with $\Cal H^+$ replaced by $\Cal H^-_{-1}$
or   $\Cal H_{-1}$ are similar, and one can also define spaces of symbol-kernels with
$\Cal H^+$ replaced by  $\Cal S_+$ (arising from inverse Fourier transformation)  etc., cf.\ \cite{G16}.

It follows from \cite{G16} Th.\ 2.6, that when $P$ is a classical $\psi $do of order $2a$,
even (cf.\ (1.2)), and elliptic avoiding a ray, then the principal
symbol $q_0(x,\xi )=p_0(x,\xi )[\xi ]^{-2a}$ of $Q=\Xi _-^{-a}P\Xi _+^{-a}$ has a factorization 
$$
q_0(x,\xi )= q_0^-(x,\xi ) q_0^+(x,\xi ),\tag A.6
$$
where 
$q_0^+=1+f^+$ with $f^+\in S_{1,0}^0({\Bbb
R}^n,{\Bbb R}^{n-1},\Cal H^+)$, and
$q_0^-=s_0+f^-$ with $f^-\in S_{1,0}^0({\Bbb R}^n,{\Bbb R}^{n-1},\Cal
H^-_{-1})$; we have here included $s_0(x)=p_0(x,0,1)$ as a factor in $q_0^-$.

The functions
$q_0^\pm$ are generalized $\psi $do symbols, like $([\xi
']\pm i\xi _n)^d$, and therefore define operators with have slightly weaker mapping
properties than $\psi $do's (see \cite{G16} Sect.\ 2), but when Poisson- and trace
operators are constructed from such symbols (see further below), they
enter in the same way as
$\psi $do symbols with the transmission property, because it is the
$\Cal H_{-1}$-property that is used. 

There is also a factorization construction where the $\Xi ^t_{\pm}$ are
replaced by the true $\psi
$do families $\Lambda _\pm^t$. Here we define
$$
Q_1=\Lambda _-^{-a}P\Lambda _+^{-a},\quad q_1(x,\xi )=\lambda
_-^{-a}\#p(x,\xi )\lambda _+^{-a},\tag A.7
$$
cf.\ (2.3). The principal symbol of $Q_1$ is again $q_0$, that factorizes as in
(A.6). By \cite{G16} Th. 2.7  there is a factorization of the full symbol
$q_1(x,\xi )$,
 leading to a factorization of $Q_1$, hence of $P$:
$$
q_1\sim q_1^-\#q_1^+,\quad Q_1\sim Q_1^-Q_1^+,
\quad P\sim \Lambda _-^{a}Q_1^-Q_1^+\Lambda _+^{a},\tag A.8
$$
where $Q_1^\pm=\OP(q_1^\pm)$ (also here, $s_0$ is included as a factor
to the left in $q_1^-$).

There is a similar full factorization in terms of $Q=\Xi _-^{-a}P\Xi _+^{-a}$.
Th.\ 2.7 of \cite{G16} does not apply directly to $\chi
_-^{-a}\#p(x,\xi )\chi _+^{-a}$ since this already passes outside the
true $\psi $do's, but the proof can be adapted to it since it is the $\Cal
H_{-1}$-properties that are used. Alternatively, we can use that $\Lambda 
_\pm^{\mu }=\Xi _\pm^{  \mu }(I+\Psi _{\pm}^{\mu })$, where $\Psi _\pm^\mu $ have
symbols $\psi _\pm^\mu(\xi ) $ in $S^0(\Cal H^\pm_{-1})$, cf.\  \cite{G15},
(1.16)ff.\ and Lemma 6.6. This  gives by insertion in (A.8), taking $\mu =a$:
$$
\aligned
&P\sim \Xi _-^a(1+\Psi _-^a)Q_1^-Q_1^+(1+\Psi _+^a)\Xi _-^{a}=\Xi
_-^aQ^-Q^+\Xi _-^{a},\text{ where}\\ 
&Q^-=(1+\Psi _-^a)Q_1^-,\quad
Q^+=Q_1^+(1+\Psi _+^a).
\endaligned\tag A.9
$$
Here $Q^-$ is of the form $s_0+F^-$ where $F^-$ has symbol in
$S^0(\Cal H^-_{-1})$, and $Q^+=1+F^+$ where   $F^+$ has symbol in
$S^0(\Cal H^+)$. 

Since the symbols $q_1^-, q_1^+, q^-, q^+$ have invertible principal
symbols, they have full parametrix symbols  $\widetilde q_1^-,\tilde
q_1^+, \widetilde q^-,\widetilde q^+$ (as in \cite{G16} Th.\ 2.9); hereby   $Q_1^-, Q_1^+, Q^-, Q^+$  have parametrices $\widetilde
Q_1^-,\widetilde  Q_1^+,\widetilde  Q^-,\widetilde Q^+$.

Pseudodifferential operators defined from a symbol $p(x,\xi )$ as in
(2.2) are a special case of $\psi $do's defined from amplitude
functions $r(x,y,\xi )$, also called symbols ``in $(x,y)$-form'', by
$$
Ru=\operatorname{OP}(r(x,y,\xi ))u 
=(2\pi )^{-n}\int_{{\Bbb R}^{2n}} e^{i(x-y)\cdot\xi
}r(x,y,\xi ) u(y)\, dyd\xi. \tag A.10
$$
Symbols $p(x,\xi )$ used as in (2.2) are then said to be ``in
 $x$-form'', and there is an asymptotic formula for the $x$-form symbol
 producing the same operator as the symbol $r(x,y,\xi )$ in (A.10) (see e.g.\
 \cite{G09} Th.\ 7.13). A useful special case is when $r(x,y,\xi )$ in (A.10) is independent of
 $x_n$ and $y'$ (which can then be left out); then it is said to be ``in
 $(x',y_n)$-form''. The transition from $x$-form to $(x',y_n)$-form goes
 as follows: An operator $R=\OP(r(x,\xi ))$ has a symbol $r'$ in
$(x',y_n)$-form related to the symbol $r$ in $x$-form by
$$
r'(x',y_n,\xi )\sim \sum_{j\in {\Bbb N}_0}\tfrac 1{j!}(i\partial_{\xi
 _n})^j\partial_{x_n}^jr(x',x_n,\xi )|_{x_n=y_n}.\tag A.11
$$

We now recall explicitly some rules from the Boutet de Monvel calculus
that we need,
cf.\ e.g.\
\cite{G09}, (10.22), Ths.\ 10.24 and 10.25: 

Trace operators arise in particular from compositions $\gamma
_0r^+Re^+$, where $R$ is a $\psi $do having the 0-transmission property.
When $R=\OP(r(x,\xi
))$, then $\gamma _0r^+Re^+$ is the trace operator with symbol $h^-r(x',0,\xi )$,
$$
\gamma _0r^+\OP(r(x,\xi ))e^+=\operatorname{OPT}(h^-r(x',0,\xi )).\tag
A.12
$$

The composition  $\varphi \mapsto r^+R(\varphi
(x')\otimes\delta (x_n))$ gives rise to a Poisson operator. It has the simplest expression when $R$ is in
$(x',y_n)$-form, namely
$$
 r^+\OP(r'(x',y_n,\xi ))(\varphi
(x')\otimes\delta (x_n))=\operatorname{OPK}(h^+r'(x',0,\xi ))\varphi
.\tag A.13
$$

These rules extend to the generalized $\psi $do's with symbols in $S^m(\Cal
H_{-1})$, $m$  integer $\le 0$, since it is the property of the symbol
being in  $\Cal H_{-1}$ with respect to $\xi _n$ that is used in the
proofs of (A.12), (A.13).

Further composition rules in the calculus are amply described in \cite{G09}
Sect.\ 10.4. Let us just recall that when $T=\operatorname{OPT}(t(x',\xi ))$ and $K=\operatorname{OPK}(k(x',\xi ))$, then
$S=TK$ is a $\psi $do $S=\OP'(s(x',\xi '))$ on ${\Bbb R}^{n-1}$ with symbol
$$
s(x',\xi ')=\tfrac1{2\pi }\int^+ t(x',\xi )\# k(x',\xi )\, d\xi
_n \tag A.14
$$
(the Leibniz product $\#$ (2.3) is used in the $(x',\xi ')$-variables).
 Here the {\it plus-integral} $\int^+f(t)\,dt$ stands for  a linear extension of, on one hand, the ordinary integral on
 ${\Bbb R}$ of functions  $f\in L_1({\Bbb R})$ and, on the other hand,  the integral over a
 curve encircling the poles in ${\Bbb C}_+$ when $f(t)$ extends to a
 meromorphic function of $t$ in ${\Bbb C}_+$. The plus-integral vanishes on
 $\Cal H^-$. %polynomials and on functions in $\Cal H^-_{-1}$. 
 It also vanishes on functions
 in $\Cal H^+$ that are $O(t^{-2})$ at infinity, since the integral can
 then be turned into an integral over a contour in ${\Bbb C}_-$. 

Calculating a plus-integral is essentially the Fourier transform of
taking a boundary value from the right: When $u(x_n)\in \Cal S(\crp)$
and $f(\xi _n)=\Cal F_{x_n\to \xi _n}(e^+u)$, then
$$
\tfrac1{2\pi }\int^+ f(\xi _n)\, d\xi
_n=\gamma _0u=\lim_{x_n\to 0+}u(x_n). \tag A.15
$$
 
Also for Poisson and trace operators, one can define them in
$x'$-form, $(x',y')$-form or $y'$-form. The adjoints of the Poisson operators $\operatorname{OPK}(k(x',\xi ))$ are the
trace operators $\operatorname{OPT}(\bar k(y',\xi ))$; the latter are
exactly the trace operators of class zero, written in $y'$-form.

We shall need the result that $\gamma _0^{a-1}u$ is described equally well
using $\Xi _+^{a-1}$ or $\Lambda _+^{a-1}$. By the remarks before (A.9),
%Recall from for any $\mu $, that $\lambda
%_+^\mu \chi  _+^{-\mu }=1+q^\mu _+$, where $q^\mu _+\in S^0(\Cal
%H^+)$. Using this for $\mu =a-1$, 
%we have that 
$\Lambda
_+^{a-1}=(1+\Psi _+^{a-1})\Xi  ^{a-1}_+$, where 
$\Psi _+^{a-1}$ has symbol $\psi _+^{a-1}(\xi )\in S^0(\Cal H^+)$. Let $u\in H^{(a-1)(s)}(\crnp)$,
$s>a-\frac12$, then $\tilde u=r^+\Xi _+^{a-1}u\in \ol H^{s-a+1}(\rnp)$
with $u=\Xi _+^{1-a}e^+\tilde u$ (by definition). We know from
\cite{G15}, Sect.\ 5, and it is reproved in Section 2 above, that
$\gamma _0^{a-1}u=\gamma _0\Xi _+^{a-1}u$ (as a boundary value from
$\rnp$). Then moreover,
$$
\gamma _0(\Lambda _+^{a-1}u)=\gamma _0((I+\Psi _+^{a-1} )\Xi _+^{a-1}u)=\gamma
_0(\Xi _+^{a-1}u)=\gamma _0^{a-1}u.\tag A.16
$$
Namely, $\gamma _0(\Psi _+^{a-1}\Xi _+^{a-1}u)=\gamma _0(\Psi _+^{a-1}e^+\tilde u)$ is in the $\xi _n$-variable
calculated as a plus-integral of $\psi ^{a-1}_+(\xi )$ multiplied by
$\widehat{e^+\tilde u}$, where both functions are in $\Cal H^+$, hence
the product is $O(\ang{\xi _n}^{-2})$ and in $\Cal H^+$, whereby the plus-integral gives
zero.
(This kind of calculation also enters in \cite{G16}, around (3.31).)

\enddocument